\documentclass[12pt]{article}

\usepackage{setspace}
\usepackage{bbm}
\usepackage{graphicx}

\usepackage{amssymb}
\usepackage{amsfonts}
\usepackage{mathrsfs}
\usepackage{amssymb}

\usepackage{amsmath}
\numberwithin{equation}{section}

\allowdisplaybreaks

\newtheorem{lemma}{Lemma}[section]

\newtheorem{theorem}{Theorem}[section]
\newtheorem{proposition}{Proposition}[section]
\newtheorem{definition}{Definition}[section]

\newtheorem{remark}{Remark}[section]
\newenvironment{proof}{{\em Proof:}}{\hfill\rule{2mm}{2mm}}

\newcommand{\dblon}{\addtolength{\baselineskip}{.1in}}

\usepackage{fancyhdr}

\begin{document}
\large
\title{ Strong convergence rate in averaging principle for stochastic hyperbolic-parabolic equations with two time-scales
\thanks{This research is funded by  National Natural
Science Foundation of China Grant (Nos. 11301403, 61573011,
11271295,
 11271013, 11301197), Science and Technology Research Projects of
Hubei Provincial Department of Education (No: D20131602) and
Foundation of Wuhan Textile University (No 2013, 2015).}}

\author{Hongbo Fu, Li Wan \\ \small College of Mathematics and Computer Science \\ \small Wuhan Textile University, Wuhan, 430073, PR
China.
\\ \small
hbfu@wtu.edu.cn (Hongbo Fu), \;wanlinju@aliyun.com (Li Wan)\\ \\ Jicheng Liu, Xianming Liu\footnote{Corresponding author}\\
\small School of Mathematics and Statistics\\\small Huazhong
University of Science and Technology, Wuhan, 430074, PR China.\\
\small jcliu@hust.edu.cn (Jicheng Liu), \;xmliu@hust.edu.cn
(Xianming Liu)} \dblon

\maketitle \small

\begin{abstract}
In this article, we investigate averaging principle for stochastic
hyperbolic-parabolic equations with two time-scales, in which both
the slow and fast components are perturbed by multiplicative noises.
Particularly, we prove that the rate of strong convergence for the
slow component to the averaged dynamics is of order $1/2$, which
significantly improves the order $1/4$ established in our previous
work.\\


\end{abstract}

\section{Introduction}
\label{intro}

This paper, which is a sequel to  Fu et al. \cite{Fu-Liu-2}, is
devoted to the strong convergence rate in averaging principle for a
coupled hyperbolic-parabolic system  with two widely separated
time-scales.

Let $D=(0, L)\subset \mathbb{R}$ be a bounded open interval. For
fixed $T_0>0$, we are concerned with the following stochastic
hyperbolic-parabolic equations with multiplicative noises:
\begin{eqnarray}
&&\!\!\!\!\!\!\frac{\partial^2 X_t^{\epsilon}(\xi)}{\partial
t^2}=\Delta X_t^{\epsilon}(\xi)+f(X_t^{\epsilon}(\xi),
Y_t^{\epsilon}(\xi))
+\sigma(X_t^{\epsilon}(\xi))\dot{W_t^1}(\xi),\label{Slow}\\
&&\!\!\!\!\!\!\frac{\partial Y_t^{\epsilon}(\xi)}{\partial
t}=\frac{1}{\epsilon}\Delta
Y_t^{\epsilon}(\xi)+\frac{1}{\epsilon}g(X_t^{\epsilon}(\xi),
Y_t^{\epsilon}(\xi))
+\frac{1}{\sqrt{\epsilon}}b(X_t^{\epsilon}(\xi),Y_t^{\epsilon}(\xi))\dot{W_t^2}(\xi),\label{Fast}\\
&&\!\!\!\!\!\!X_t^{\epsilon}(0)=Y_t^{\epsilon}(0)=X_t^{\epsilon}(l)=Y_t^{\epsilon}(l)=0, \label{Boundary}\\
&&\!\!\!\!\!\!X_0^{\epsilon}(\xi)=X_0(\xi),
Y_0^{\epsilon}(\xi)=Y_0(\xi), \frac{\partial
X_t^{\epsilon}(\xi)}{\partial t}\Big|_{t=0}=\dot{X}_0(\xi),
\label{Initial}
\end{eqnarray}
where the space variable $\xi\in D$, the time $t\in[0, T_0]$. Here
the forcing noises $W^1_t(\xi)$ and $W^2_t(\xi)$ are mutually
independent Wiener processes on a complete stochastic basis
$(\Omega, \mathscr{F}, \mathscr{F}_t, \mathbb{P})$, which will be
specified later. Also, the precise conditions on $f, g, \sigma$ and
$b$ will be presented in the next section. The positive and small
parameter $\epsilon$ measures the  ratio of the time scales  between
slow component $X_t^\epsilon(\xi)$ and fast component
$Y_t^\epsilon(\xi)$.

The system in form of \eqref{Slow}-\eqref{Initial} is an abstract
model for a random vibration of an elastic string with
 external force on a large time scale. More generally, the nonlinear coupled wave-heat
equations with fast and slow time scales may describe a
thermoelastic wave propagation in a random medium (Chow
\cite{Chow-1}), the interactions of fluid motion with other forms of
waves (Leung \cite{Leung},  Zhang and Zuazua \cite{X Zhang}), wave
phenomena which are  temperature related (Leung \cite{Leung-0}),
magneto-elasticity processes (Rivera and Racke \cite{Rivera}) as
well as biological problems (Choi and Miller \cite{Choi}, Cardetti
and Choi \cite{Cardetti}  and Wu et al. \cite{S.H. Wu}).

We are often interested in the dynamical evolution of the slow
component $X^\epsilon_t(\xi)$ as the scale parameter $\epsilon$ goes
to zero. Due to the  separation of time-scale, a simplified
equation, which excludes the fast component and approximates the
dynamics of slow component, is highly desirable. Such a simplified
equation can be obtained by the so-called averaging procedure. This
result has been proved in our previous paper \cite{Fu-Liu-2} in the
case when the driving noises are additive type.

Averaging principle is a powerful tool to analyze the asymptotic
behavior for slow-fast dynamical systems. It was first formulated by
Bogoliubov \cite{Bogoliubov} for deterministic differential
equations. For its validity to stochastic differential equations
(SDEs) with Gaussian noise, we mainly refer to the well known paper
Khasminskii \cite{Khas}, the works of Freidlin and Wentzell
\cite{Freidlin-Wentzell1,Freidlin-Wentzell2}, Veretennikov
\cite{Vere1,Vere2} and Kifer \cite{Kifer1,Kifer2,Kifer3}. Further
progresses on averaging for SDEs with non-Gaussian noise have been
made in Xu and his co-workers \cite{Xu,Xu1}. Recently, there are
increasing interests to extend the classical averaging about SDEs to
the case of stochastic partial differential equations (SPDEs). The
basic results are due to Cerrai and Freidlin \cite{Cerrai1} and
Cerrai \cite{Cerrai2,Cerrai-Siam}, {in which the} averaging
principle of stochastic reaction-diffusion equation perturbed by
Gaussian noise, with no explicit convergence rate being given, has
been systematically studied. In Wang and Roberts \cite{wangwei}, the
rate of convergence in probability of the original equation to the
averaged equation has been determined. In Bao et al. \cite{Bao}, the
averaging dynamics for two-time-scale SPDEs with $\alpha-$stable
noises has been derived.

Once the averaging principle is established, an important question
arises as to  {how fast will the original slow component } converge
to the effective dynamics. In Br\'{e}hier \cite{Brehier}, when
additive noise is included only in the fast component, explicit
strong convergence rate of order $\frac{1}{2}-\varepsilon$ for
arbitrary small $\varepsilon>0$ for averaging of stochastic
parabolic equations is obtained. However, the order will be
decreased to $\frac{1}{5}$ if the noise also acts on the slow
variable  directly (see \cite[Section 1]{Brehier}). These
convergence rates can be compared with order $\frac{1}{2}$ obtained
for the finite dimensional stochastic {dynamical systems} (see,
e.g., Liu \cite{LiuDi},  Khasminskii \cite{Khas2}). An interesting
question thus occurs as to whether it is possible to get an order
$\frac{1}{2}$ for strong convergence in averaging of stochastic
dynamical systems in infinite dimensional space. In this article, we
present a positive answer to it by dealing with  the coupled
hyperbolic-parabolic equations with multiplicative noises in the
form of \eqref{Slow}-\eqref{Initial}. To be more precise, we will
show the slow component $X_t^\epsilon(\xi)$ can be approximated by
the solution of an averaged system, which enjoys strong convergence
order $\frac{1}{2}$ and is governed by a stochastic wave equation
constructed by averaging the slow dynamic with respect to the
stationary measure associated with the fast component. To the
author's knowledge, this is the first paper to obtain the strong
convergence order $\frac{1}{2}$ for averaging of SPDEs in the case
when the slow equation is subjected with noise, which obviously
exceeds order $\frac{1}{4}$ established in our previous work
\cite{Fu-Liu-2}.

To construct the averaged system, a key point is to show the
existence of an invariant measure with exponentially mixing property
for the fast equation and this can be obtained by the same
discussion as in \cite{Fu-Liu-2}, where a dissipative condition is
needed. To provide an explicit error bound on the difference between
the solution of the original system and the solution of the averaged
equation, we follow the general lines of the arguments introduced in
Liu \cite{LiuDi}. But in our case, this procedure will be a little
more technical as it involves a system   in infinite dimensional
space. In order to obtain a sharp rate for strong convergence, we
need to require  the differentiability with respect to the parameter
of the solution to fast motion equation with a frozen slow component
(see equations \eqref{Frozen-Fast}-\eqref{Frozen-Initial} in Section
\ref{section-ergodicty}). Therefore, in our setup we introduce
additional derivable conditions on drift and diffusion coefficients.

The organization of this paper is as follows: In Section 2, we
recall some basic concepts and results for later use. In Section 3,
we prove the existence, uniqueness and the energy identity for an
abstract hyperbolic-parabolic equation. In Section 4, we study the
ergodicity property for the fast component of the system
\eqref{Slow}-\eqref{Initial}. In Section 5, some priori estimates is
presented. Section 6 contains the main results of the paper as
presented in Theorem \ref{Theorem}. Finally, a necessary lemma is
proved in the last section.

\section{Preliminary}
To transform the system \eqref{Slow}-\eqref{Initial} as the abstract
evolution equations, we present some notations and recall some
well-known facts for later use.

Let $(U,(\cdot,\cdot)_U,\|\cdot \|_U )$ and
$(V,(\cdot,\cdot)_V,\|\cdot \|_V)$ be two separable Hilbert spaces.
We denote the space of bounded linear operators from $U$ to $V$ by
$\mathcal{L}(U, V)$ with the usual operator norm
$\|\cdot\|_{\mathcal{L}(U, V)}$. The space of Hilbert-Schmidt
operators from $U$ to $V$ is denoted by $\mathcal{L}_2(U, V)$, which
is equipped with norm
\[\|T\|_{\mathcal{L}_2(U, V)}=\left(\sum\limits_{k\in \mathbb{N}}\|Tu_k\|^2_V\right)^{1/2},\]
where $\{u_k\}_{k\in\mathbb{N}}$ is an arbitrary orthonormal basis
of $U$. For simplicity, we write $\mathcal{L}(U)=\mathcal{L}(U, U)$,
and $\mathcal{L}_2(U)=\mathcal{L}_2(U, U)$.  Let $Q\in
\mathcal{L}(U)$ be a nonnegative and symmetric operator and use
$Q^{\frac{1}{2}}$ to denote   the unique positive square root of
$Q$. The space of Hilbert-Schmidt operators from $Q^{\frac{1}{2}}U$
to $U$ is denoted by $\mathcal{L}_{2, Q}^0(U)$, equipped with the
scalar product \[\Big(S, T
\Big)_Q=\sum\limits_{k\in\mathbb{N}}\big(SQ^\frac{1}{2}u_k,TQ^\frac{1}{2}u_k\big)_U.\]
The associated  norm is given by
$$\|T\|_{Q}=\|TQ^\frac{1}{2}\|_{\mathcal{L}_2(U)}.$$
For a fixed domain $D=(0, L)$, let $H$ be the Hilbert space
$L^2(D)$, endowed with the usual scalar product $(\cdot, \cdot)_H$
and the corresponding norm $\|\cdot\|$.

Let $\{ e_k(\xi)\}_{k\in\mathbb{N}}$ denote the complete
orthornormal system of eigenfunctions in $H$ such that, for $k =
1,2,\ldots$,
\begin{equation*}\label{eigenfunction} -\Delta
e_k(\xi)=\alpha_ke_k(\xi),\;\;e_k(0)=e_k(L)=0 ,
\end{equation*}
with $0<\alpha_1\leq\alpha_2\leq\cdots\alpha_k\leq\cdots$. {Here we
would like to recall the fact that
$e_k(\xi)=\sqrt{\frac{2}{L}}\sin\frac{k\pi\xi}{L}$ and
$\alpha_k=-\frac{k^2\pi^2}{L^2}$ for $k = 1,2,\cdots$.}

Let $A$ be the  Laplacian  operator  with $\mathcal
{D}(A)=H^1_0(D)\cap H^2(D)$ . For $s\in \mathbb{R}$, we introduce
Hilbert space $\mathcal {H}^s :=\mathcal {D}((-A)^{s/2})$, which is
equipped with the scalar product
\begin{equation*}
\langle u,v \rangle_s=\sum\limits_{k=1}^\infty
\alpha_k^s(u,e_k)_H(v,e_k)_H
\end{equation*}
and norm
\begin{equation*}
\|u\|_s=\left\{\sum\limits_{k=1}^\infty\alpha_k^s\big(u,
e_k\big)_H^2\right\}^\frac{1}{2}
\end{equation*}
for $u,v\in \mathcal {H}^s .$ It is known that $\mathcal {H}^0=H,
\mathcal {H}^1=H_0^1(D)$ and $ \mathcal {H}^2=H_0^1(D)\cap H^2(D)$
with equivalent norms. We note that in the case  of $s>0$, $\mathcal
{H}^{-s}$ can be identified with the dual space $(\mathcal
{H}^{s})^*$, i.e., the space of the  linear functionals on $\mathcal
{H}^{s}$, which are continuous with respect to the topology induced
by the norm $\|\cdot\|_s$.

Recall that the Green function $G(\xi,\zeta;t)$ for linear equation
$(\partial/\partial t-A)X(t, \xi)=0$ can be expressed as
$$G(\xi,\zeta,t)=\sum\limits_{k=1}^\infty
e^{-\alpha_kt}e_k(\xi)e_k(\zeta).$$ The associated Green's operator
is defined by, for any $h(\xi)\in H$,
\begin{equation*}
G_th(\xi)=\int_DG(\xi,\zeta,t)h(\zeta)d\zeta=\sum\limits_{k=1}^\infty
e^{-\alpha_kt}e_k(\xi)\big(e_k, h\big)_H.
\end{equation*}
It is straightforward to check that $\{G_t\}_{t\geq0}$ is a
contractive semigroup  on $H$. For the  linear wave equation
$(\partial^2/\partial t^2-A)Y(t, \xi)=0$, its Green's function is
given by
$$S(\xi,\zeta, t)=\sum\limits_{k=1}^\infty\frac{\sin\{\sqrt{\alpha_k}t\}}{\sqrt{\alpha_k}}e_k(\xi)e_k(\zeta).$$
It is easy to shown that the above series converges in $L^2(D\times
D)$ and the associated Green's operator is defined by, for any
$h(\xi)\in H$,
$$S_th(\xi)=\int_DS(\xi,\zeta, t)h(\zeta)d\zeta=\sum\limits_{k=1}^\infty
\frac{\sin\{\sqrt{\alpha_k}t\}}{\sqrt{\alpha_k}}e_k(\xi)\big(e_k,
h\big)_H.$$ For Green's operator $S_t$, we present
the following results (see Chow \cite{Chow}): 
\begin{proposition}\label{lemma-2}
Let $\phi_t(\cdot)\in L^2(\Omega\times[0, T]; H)$ be an
$\mathscr{F}_t-$adapted process which satisfies
$\mathbb{E}\int_0^T\|\phi_t(\cdot)\|^2dt<\infty.$ Then
$\nu_t(\cdot)=\int_0^tS_{t-s}\phi_s(\cdot)ds$ is a continuous,
$\mathscr{F}_t-$adapted $\mathcal {H}^1-$valued process and its time
derivative $\dot{\nu}_t(\cdot)=\frac{\partial}{\partial
t}{\nu}_t(\cdot)$ is a continuous $H-$valued process such that
\begin{eqnarray}
\mathbb{E}\sup\limits_{0\leq t\leq T}\|\nu_t(\cdot)\|^2_1\leq
T\mathbb{E}\int_0^T\|\phi_s(\cdot)\|^2ds  \label{S-t-2}
\end{eqnarray}
and
\begin{eqnarray}
\mathbb{E}\sup\limits_{0\leq t\leq T}\|\dot{\nu}_t(\cdot)\|^2\leq
T\mathbb{E}\int_0^T\|\phi_s(\cdot)\|^2ds.\label{S-t-3}
\end{eqnarray}
\end{proposition}

Next we recall the definition of   Wiener process in infinite
dimensional space. For more details, see Chow \cite{Chow} or Da
Prato and Zabczyk \cite{Daprato}. Let $W_t^1$ and $W_t^2$ be two
mutually independent $H-$valued Wiener process on a complete
stochastic basis  $(\Omega, \mathscr {F},\mathscr {F}_t,\mathbb{P})$
with  nonnegative and symmetric covariance operators $Q_1$ and
$Q_2$, respectively, which is defined by $Q_ie_k=\lambda_{i, k}e_k$
for $i=1,2$, where $\{\lambda_{i,k}\}_{k\in \mathbb{N}}$ is a
sequence of nonnegative real numbers satisfying
\begin{eqnarray}
TrQ_i=\sum\limits_{k\in\mathbb{N}}\lambda_{i, k}:=r_i<
{+\infty},\label{covariance-bound}
\end{eqnarray}
and $\{e_k\}_{k\in\mathbb{N}}$ is  the complete orthonormal system
for  $H$,  consisting  of eigenfunctions of Laplacian operator $A$.
Formally, for every $i=1,2$, the Wiener processes $W^i_t$ is defined
as the  {infinite sums}
\[W^{i}_t =\sum\limits_{k\in\mathbb{N}} {\lambda^\frac{1}{2}_{i, k}}\beta_{i, k}(t)e_k,\;t\geq 0, \]
where $\{\beta_{i, k}(t)\}^{i=1, 2}_{k\in\mathbb{N}}$ are mutually
independent real-valued Brownian motions on the same  stochastic
basis $(\Omega, \mathscr{F}, \mathscr{F}_t, \mathbb{P})$. With
condition \eqref{covariance-bound}, we know that $Q_1$ and $Q_2$ are
Hilbert-Schmidt operators. Therefore, for $i=1,2$, $W_t^i$ is often
called the ${Q_i}-$Wiener process on $H$. For simplicity of
presentation, {we will write} $\mathcal{L}^0_{2,Q_i}
={\mathcal{L}_2(Q_i^{\frac{1}{2}}(H), H)}$ and
$\|\cdot\|_{Q_i}=\|\cdot\|_{\mathcal{L}_2(Q_i^{\frac{1}{2}}(H),
H)}$.

\begin{proposition}\label{Pro-2.3}
Let $\phi_t(\cdot)\in L^2(\Omega\times[0, T]; \mathcal{L}^0_{2,Q_1}
)$ be an $\mathscr{F}_t-$adapted process which satisfies
$\mathbb{E}\int_0^T\|\phi_t(\cdot)\|^2_{Q_1}dt<\infty.$ Then
\begin{eqnarray*}
\mu_t(\cdot)=\int_0^tS_{t-s}\phi_s(\cdot)dW_s^1,\;t\in[0,T]
\end{eqnarray*}
is a continuous adapted $\mathcal {H}^1-$valued process, and its
derivative $\dot{\mu}_t(\cdot)=\frac{\partial}{\partial
t}\mu_t(\cdot)$ is a continuous process in $H$ (see  Section 5.2 in
Chow \cite{Chow} for details). Moreover, the following inequalities
hold:
\begin{eqnarray}
\mathbb{E}\|\mu_t(\cdot)\|_1^2\leq
\mathbb{E}\int_0^t\|\phi_s(\cdot)\|^2_{Q_1}ds,\label{S-t-1-stocha}
\end{eqnarray}
\begin{eqnarray}
\mathbb{E}\|\dot{\mu}_t(\cdot)\|^2\leq
\mathbb{E}\int_0^t\|\phi_s(\cdot)\|^2_{Q_1}ds.\label{S-t-2-stocha}
\end{eqnarray}
\end{proposition}

A natural way to give a rigorous meaning to system
\eqref{Slow}-\eqref{Initial} is in terms of the following integral
equations:
\begin{eqnarray}
&&X^\epsilon_t=S_t'X_0+S_t\dot{X}_0+\int_0^tS_{t-s}f(X_s^{\epsilon},
Y^{\epsilon}_s)ds+\int_0^tS_{t-s}\sigma(X_s^{\epsilon})dW^1_s,\label{Integral-Slow}\\
&&Y^\epsilon_t=G_{t/\epsilon}Y_0+\frac{1}{\epsilon}\int_0^tG_{(t-s)/\epsilon}g(X^\epsilon_s,
Y_s^\epsilon)ds\nonumber\\
&&\quad\quad\quad\quad\quad\;\;\;+\frac{1}{\sqrt{\epsilon}}\int_0^tG_{(t-s)/\epsilon}b(X^\epsilon_s,
Y_s^\epsilon)dW^2_s,\label{Integral-Fast}
\end{eqnarray}
where $S'_t=\frac{d}{dt}S_t$ is the derived Green's operator with
integral kernel
$$K'(\xi,\zeta, t)=\sum\limits_{k=1}^\infty\cos\{\sqrt{\alpha_k}t\}e_k(\xi)e_k(\zeta).$$
As a solution to system \eqref{Slow}-\eqref{Initial}, we take the
so-called mild sense.
\begin{definition}
If  $(X^\epsilon_t, Y^\epsilon_t)$ is an adapted process over
$(\Omega, \mathscr{F}, \mathscr{F}_t, \mathbb{P})$ such that
$\mathbb{P}-a.s.$ the
 integral equations \eqref{Integral-Slow}-\eqref{Integral-Fast}
 hold true for all $t>0$, we say
that it is a {mild solution} for equations
\eqref{Slow}-\eqref{Initial}.
\end{definition}

Let us introduce the following set of additional assumptions.

(A1) For the mapping $f:H\times H\rightarrow H$, we require that
there exists a constant $L_f$ such that for any $x, y,h,k\in H$ the
following directional derivatives are well-defined and satisfy
\begin{eqnarray}
&&\|D_xf(x,y)\cdot h\|\leq L_f\|h\|,\label{f-1}\\
&&\|D_yf(x,y)\cdot h\|\leq L_f\|h\|,\label{f-2}\\
&&\|D^2_{xy}f(x,y)\cdot (h,k)\|\leq L_f\|h\|\cdot\|k\|\label{f-3},\\
&&\|D^2_{yy}f(x,y)\cdot (h,k)\|\leq L_f\|h\|\cdot\|k\|.\label{f-4}
\end{eqnarray}
Moreover, we assume that $f$ is bounded.

(A2) Assume the mapping $\sigma: H\rightarrow \mathcal{L}_{2, Q_1}^0
$ satisfies the global Lipschitz condition and the sublinear growth.

(A3) For the mapping $g:H\times H\rightarrow H$, we assume that
there exist two constants $C_g$ and $L_g$ such that for any $x,
y,h,k\in H$ the following directional derivatives are well-defined
and satisfy
\begin{eqnarray}
&&\|D_xg(x,y)\cdot h\|\leq C_g\|h\|,\label{g-1}\\
&&\|D_yg(x,y)\cdot h\|\leq L_g\|h\|.\label{g-2}
\end{eqnarray}

(A4) For  the mapping $b: H\times H\rightarrow \mathcal{L}_{2,
Q_2}^0$, we require that there exist two constants $C_b$ and $L_b$
such that for any $x, y,h \in H$ the following directional
derivatives are well-defined and satisfy
\begin{eqnarray}
&&\|D_xb(x,y)\cdot h\|_{Q_2}\leq C_b\|h\|,\label{b-1}\\
&&\|D_yb(x,y)\cdot h\|_{Q_2}\leq L_b\|h\|.\label{b-2}
\end{eqnarray}

(A5) Assume the fast motion equation satisfies the strong
dissipative condition, that is
\begin{equation}\label{decay}
{\kappa:=2\alpha_1-2L_g- L^2_{b}>0. }
\end{equation}

\begin{remark}
A simple example of the drift coefficient  $f$ can be given by
\begin{equation*}
f(x,y)=f_1(x)+f_2(y),
\end{equation*}
here $f_1, f_2: H\rightarrow H$ are of class $C^2$, bounded and with
bounded derivatives up to the second order.
\end{remark}

\begin{remark}
{With assumptions} (A1) and (A3), it is immediate to check that $f$
and $g$  are Lipschitz mappings from $H\times H $ to $H$. In the
same way, by assumption (A4), the mapping $b: H\times H\rightarrow
\mathcal{L}_{2, Q_2}^0$ is Lipschitz-continuous.
\end{remark}

\begin{remark}
For the convenience of notations, the symbols $C$ with or without
subscripts will denote a positive constant that is unimportant and
may have different values from one line to another one in the
sequel.
\end{remark}

\section{Existence, uniqueness and  energy equality}
In this section, we prove the existence and uniqueness of mild
solutions to Eqs. \eqref{Integral-Slow}-\eqref{Integral-Fast}.   Let
$U_t$ and $V_t$ be two $\mathscr{F}_t-$adapted  processes in $H$,
and let $\tilde{U}_t$ and $\tilde{V}_t$ be two
$\mathscr{F}_t-$adapted processes in $\mathcal{L}_{2, Q_1}^0$ and
$\mathcal{L}_{2, Q_2}^0$, respectively, such that
\begin{equation}\label{linear-integrable}
\mathbb{E}\int_0^{T_0}(\|U_s\|^2+\|\tilde{U}_s\|^2_{Q_1}+\|V_s\|^2+\|\tilde{V}_s\|^2_{Q_2})ds<
\infty.
\end{equation}
Now, for fixed $x_0\in \mathcal {H}^1, \dot{x}_0, y_0\in H$ we first
consider the integrals
\begin{eqnarray}
&&X_t=S_t'x_0+S_t\dot{x}_0+\int_0^tS_{t-s}U_sds+\int_0^tS_{t-s}\tilde{U}_sdW^1_s\label{linear-Integral-Slow}
\end{eqnarray}
and
\begin{eqnarray}
&&Y_t=G_{t}y_0+\int_0^tG_{t-s}V_sds+
\int_0^tG_{t-s}\tilde{V}_sdW^2_s.\label{linear-Integral-Fast}
\end{eqnarray}
They are respectively mild solutions to the linear equations with
additive noise
\begin{eqnarray}
\frac{d^2}{dt^2}X_t=AX_t+U_t+\tilde{U}_t\dot{W}_t^1,
X_0=x_0,\frac{dX_t}{dt}\Big|_{t=0}=\dot{x}_0\nonumber
\end{eqnarray}
and
\begin{eqnarray}
\frac{d}{dt}Y_t=AY_t+V_t+\tilde{V}_t\dot{W}_t^2, Y_0=y_0.\nonumber
\end{eqnarray}
By Chow \cite[Theorem 3.5, Chapter 5]{Chow}, for any $T_0>0$ the
linear problem \eqref{linear-Integral-Slow} has a unique solution
$$X_t\in L^2(\Omega; {C}([0, T_0]; \mathcal {H}^1))$$ with time
derivatve $\dot{X}_t=\frac{d}{dt}X_t\in L^2(\Omega; {C}([0, T_0];
H))$. Moreover, we have  the energy equality
\begin{eqnarray}
\|\dot{X}_t\|^2+\|X_t\|^2_1&=&\|\dot{x}_0\|^2
+\|x_0\|_1^2+2\int_0^t\big(\dot{X}_s,
U_s\big)_Hds+2\int_0^t\big(X_s, \tilde{U}_sdW_s^1\big)_H\nonumber\\
&&+\int_0^t\|\tilde{U}_s\|^2_{Q_1}ds,\;\;
a.s..\label{energy-linear-slow}
\end{eqnarray}
By Chow \cite[Theorem 5.3, Chapter 3]{Chow}, the linear problem
\eqref{linear-Integral-Fast} has a unique solution which is a
process in $\mathcal {H}^1$ with continuous path in $H$ such that
the energy equality holds true:
\begin{eqnarray}\label{energy-linear-fast}
\|Y_t\|^2 &=&\|y_0\|^2+2\int_0^t\langle AY_s, Y_s\rangle
ds+2\int_0^t\big(Y_s, V_s\big)_Hds+2\int_0^t\big(Y_s,
\tilde{V}_sdW^2_s\big)_H\nonumber\\
&&+\int_0^t\|\tilde{V}_s\|^2_{Q_2}ds,\;\; a.s.,
\end{eqnarray}
where $\langle , \rangle$ denotes the dualization between
$\mathcal{H}^{-1}$ and $\mathcal{H}^1$. Now, for fixed $X_0\in
\mathcal {H}^1, \dot{X}_0, Y_0\in H$ we consider the nonlinear
problems
\begin{eqnarray}
&&X_t=S_t'X_0+S_t\dot{X}_0+\int_0^tS_{t-s}f(X_s, Y_s)ds+\int_0^tS_{t-s}\sigma(X_s)dW^1_s,\label{Nonlinear-Integral-Slow}\\
&&Y_t=G_{t}Y_0+\int_0^tG_{t-s}g(X_s, Y_s)ds+ \int_0^tG_{t-s}{b(X_s,
Y_s)}dW^2_s.\label{Nonlinear-Integral-Fast}
\end{eqnarray}
\begin{theorem}
{Suppose} that assumptions (A1)-(A4) are satisfied. For any $X_0\in
\mathcal {H}^1, \dot{X}_0, Y_0\in H$ and $T_0>0$ the system
\eqref{Nonlinear-Integral-Slow}-\eqref{Nonlinear-Integral-Fast} has
a unique mild solution $$(X_t, Y_t)\in L^2(\Omega; {C}([0, T_0];
\mathcal {H}^1))\times L^2(\Omega; {C}([0, T_0]; H))$$ with
$\dot{X}_t =\frac{d}{dt}X_t\in L^2(\Omega; {C}([0, T_0]; H))$.
Moreover, the solution enjoys the energy equalities
\begin{eqnarray}
\!\!\!\!\!\!\|\dot{X}_t\|^2+\|X_t\|^2_1&=&\|\dot{X}_0\|^2+\|X_0\|_1^2+2\int_0^t\big(\dot{X}_s,
f(X_s, Y_s)\big)_Hds \nonumber\\
&&\!\!\!\!\!\!\!\!\!+2\int_0^t\big(X_s,
\sigma(X_s)dW_s^1\big)_H+\int_0^t\|\sigma(X_s)\|^2_{Q_1}ds,\;\;a.s.
\label{energy-nonlinear-slow}
\end{eqnarray}
and
\begin{eqnarray}\label{energy-nonlinear-fast}
\!\!\!\!\!\!\|Y_t\|^2 &=&\|Y_0\|^2+2\int_0^t\langle AY_s, Y_s\rangle
ds+2\int_0^t\big(Y_s, g(X_s,
Y_s)\big)_Hds\nonumber\\
&+& 2\int_0^t\big(Y_s,b(X_s, Y_s) dW^2_s\big)_H+\int_0^t\|b(X_s,
Y_s)\|^2_{Q_2}ds,\;\;a.s..
\end{eqnarray}
\begin{proof}
 We will verify the existence by a standard approximation
argument. Let
\begin{eqnarray*}
&&X_t^0=X_0+\int_0^t\dot{X}_0ds,\\
&&Y_t^0=Y_0, t\geq 0.
\end{eqnarray*}
For $n\in \mathbb{N},$ let $(X^n_t, Y^n_t)$ be the unique solution
to the linear equations
\begin{eqnarray}
&&X_t^n=S_t'X_0+S_t\dot{X}_0+\int_0^tS_{t-s}f(X_s^{n-1}, Y_s^{n-1})ds\nonumber\\
&&\qquad\;+\int_0^tS_{t-s}\sigma(X_s^{n-1}, Y_s^{n-1})dW^1_s,\label{Nonlinear-Appro-Slow}\\
&&Y_t^n=G_{t}Y_0+\int_0^tG_{t-s}g(X_s^{n-1},
Y_s^{n-1})ds\nonumber\\
&&\qquad+ \int_0^tG_{t-s}b(X_s^{n-1},
Y_s^{n-1})dW^2_s.\label{Nonlinear-Appro-Fast}
\end{eqnarray}
Note that, for any $n\in \mathbb{N}$, the existence of  $(X^n_t,
Y^n_t)$ follows from the existence of linear problems
\eqref{linear-Integral-Slow} and \eqref{linear-Integral-Fast}.
Additional, the time derivative of ${X}_t^n$ satisfies
\begin{eqnarray}
&&\dot{X}_t^n=S_t''X_0+S_t'\dot{X}_0+\int_0^tS'_{t-s}f(X_s^{n-1}, Y_s^{n-1})ds\nonumber\\
&&\qquad\;+\int_0^tS'_{t-s}\sigma(X_s^{n-1},
Y_s^{n-1})dW^1_s,\nonumber
\end{eqnarray}
where $S''_t$ denotes the second order derivative of Green's
operator with integral kernel
$$K''(\xi,\zeta, t)=-\sum\limits_{k=1}^\infty\sqrt{\alpha_k}\sin\{\sqrt{\alpha_k}t\}e_k(\xi)e_k(\zeta).$$
Our aim is  to show that $\{(X^n_t, Y^n_t)\}_{n\in\mathbb{N}}$ forms
a Cauchy sequence in $ L^2(\Omega; {C}([0, T_0]; \mathcal
{H}^1))\times L^2(\Omega; {C}([0, T_0]; H))$. For any $t\in[0,
T_0]$, the energy equality \eqref{energy-linear-slow} yields
\begin{eqnarray}
&&\|\dot{X}_t^{n+1}-\dot{X}_t^{n}\|^2+\|X_t^{n+1}-X_t^{n}\|_1^2\nonumber\\
&&=2\int_0^t\big(\dot{X}_s^{n+1}-\dot{X}_s^{n}, f(X_s^n,
Y^n_s)-f(X_s^{n-1}, Y^{n-1}_s)\big)_Hds\nonumber\\
&&+2\int_0^t\big({X}_s^{n+1}-{X}_s^{n},
(\sigma(X_s^n)-b(X_s^{n-1}))dW_s^1\big)_H\nonumber\\
&&+\int_0^t\|\sigma(X_s^n)-\sigma(X_s^{n-1})\|^2_{Q_1}ds.\nonumber
\end{eqnarray}
This equality, together with the B-D-G inequality, allows us to get
\begin{eqnarray}
&&\mathbb{E}\sup\limits_{0\leq s\leq
t}\left[\|\dot{X}_s^{n+1}-\dot{X}_s^{n}\|^2+\|X_s^{n+1}-X_s^{n}\|_1^2\right]\nonumber\\
&&\leq
C\mathbb{E}\Big[\int_0^t\|\dot{X}_s^{n+1}-\dot{X}_s^{n}\|^2ds+\int_0^t\|{X}_s^{n+1}-{X}_s^{n}\|^2ds\nonumber\\
&&\quad+\int_0^t\|{X}_s^{n}-{X}_s^{n-1}\|^2ds+\int_0^t\|{Y}_s^{n}-{Y}_s^{n-1}\|^2ds\Big].\label{Energy-Slow-1}
\end{eqnarray}
Note that for $s\in [0, t]$, we have
\begin{eqnarray*}
\dot{X}_s^{n+1}-\dot{X}_s^n&=&\int_0^sS'_{s-r}\left(f(X_r^{n},Y_r^{n})-f(X_r^{n-1},Y_r^{n-1})\right)dr\\
&+&\int_0^sS'_{s-r}\left(\sigma(X_r^{n})-\sigma(X_r^{n-1})\right)dW_r^1,
\end{eqnarray*}
and then, due to the  contractive property of $\{S'_t\}_{t\geq 0}$
on $H$, by H\"{o}lder's inequality we obtain
\begin{eqnarray}
\mathbb{E}\|\dot{X}_s^{n+1}-\dot{X}_s^n\|^2&\leq&
s\int_0^s\mathbb{E}\|f(X_r^{n},Y_r^{n})-f(X_r^{n-1},Y_r^{n-1})\|^2dr\nonumber\\
&+&C\int_0^s\mathbb{E}\|\sigma(X_r^{n})-\sigma(X_r^{n-1})\|^2dr\nonumber\\
&\leq&C_{T_0}\int_0^s\|X_r^{n}-X_r^{n-1}\|^2+\|Y_r^{n}-Y_r^{n-1}\|^2dr.\nonumber
\end{eqnarray}
By integrating over the interval $[0,t]$ and a change of the order
of integration, we get
\begin{eqnarray}
\mathbb{E}\int_0^t\|\dot{X}_s^{n+1}-\dot{X}_s^n\|^2ds&\leq&
C_{T_0}\mathbb{E}\int_0^t\int_0^s(\|X_r^{n}-X_r^{n-1}\|^2+\mathbb{E}\|Y_r^{n}-Y_r^{n-1}\|^2)drds\nonumber\\
&=&C_{T_0}\mathbb{E}\int_0^t\int_r^t(\|X_r^{n}-X_r^{n-1}\|^2+\mathbb{E}\|Y_r^{n}-Y_r^{n-1}\|^2)dsdr\nonumber\\
&=&C_{T_0}\mathbb{E}\int_0^t(t-r)
(\|X_r^{n}-X_r^{n-1}\|^2+\mathbb{E}\|Y_r^{n}-Y_r^{n-1}\|^2)dr\nonumber\\
&\leq&C_{T_0}\mathbb{E}\int_0^t
(\|X_r^{n}-X_r^{n-1}\|^2+\mathbb{E}\|Y_r^{n}-Y_r^{n-1}\|^2)dr.\nonumber
\end{eqnarray}
Now, due to  \eqref{Energy-Slow-1}, this yields
\begin{eqnarray}
&&\mathbb{E}\sup\limits_{0\leq s\leq
t}\left[\|\dot{X}_s^{n+1}-\dot{X}_s^{n}\|^2+\|X_s^{n+1}-X_s^{n}\|_1^2\right]\nonumber\\
&&\leq C\mathbb{E}\int_0^t\big[\|{X}_s^{n+1}-{X}_s^{n}\|^2
+\|{X}_s^{n}-{X}_s^{n-1}\|^2\nonumber\\
&&\qquad+\|{Y}_s^{n}-{Y}_s^{n-1}\|^2\big]ds
\nonumber\\
&&\leq C\mathbb{E}\int_0^t\big[\|{X}_s^{n+1}-{X}_s^{n}\|^2_1
+\|{X}_s^{n}-{X}_s^{n-1}\|^2_1\nonumber\\
&&\qquad+\|{Y}_s^{n}-{Y}_s^{n-1}\|^2\big]ds\label{Energy-Slow-1-1}
\end{eqnarray}
By energy equality \eqref{energy-linear-fast}, assumption (A3) and
the fact $\langle
A(Y_s^{n+1}-Y_s^{n}),Y_s^{n+1}-Y_s^{n}\rangle\leq0$, we obtain
\begin{eqnarray}
\|Y_t^{n+1}-Y_t^{n}\|^2&\leq&
2\mathbb{E}\int_0^t\big(Y_s^{n+1}-Y_s^{n}, g(X_s^n,
Y^n_s)-g(X_s^{n-1},
Y^{n-1}_s)\big)_Hds\nonumber\\
&+&2\mathbb{E}\int_0^t\big(Y_s^{n+1}-Y_s^{n}, \left(b(X_s^n,
Y^n_s)-b(X_s^{n-1},
Y^{n-1}_s)\right)dW_s^2\big)_H\nonumber\\
&+&\mathbb{E}\int_0^t\|b(X_s^n,Y_s^n)-b(X_s^{n-1},Y_s^{n-1})\|_{Q_2}^2ds.\nonumber
\end{eqnarray}
By B-D-G inequality, we have
\begin{eqnarray}
\mathbb{E}\sup\limits_{0\leq s\leq t}\|Y_s^{n+1}-Y_s^{n}\|^2&\leq&
C\mathbb{E}\int_0^t\left[\|Y_s^{n+1}-Y_s^n\|^2+\|X_s^{n}-X_s^{n-1}\|^2+\|Y_s^{n}-Y_s^{n-1}\|^2\right]ds\nonumber\\
&\leq&
C\mathbb{E}\int_0^t\left[\|Y_s^{n+1}-Y_s^n\|^2+\|X_s^{n}-X_s^{n-1}\|^2_1+\|Y_s^{n}-Y_s^{n-1}\|^2\right]ds.\nonumber
\end{eqnarray}
 According to \eqref{Energy-Slow-1-1}, we get
\begin{eqnarray}
&&\mathbb{E}\sup\limits_{0\leq r\leq
t}\left[\|X_r^{n+1}-X_r^{n}\|_1^2+\|Y_r^{n+1}-Y_r^{n}\|^2\right]\nonumber\\
&&\leq
C\mathbb{E}\left[\int_0^t(\|X_s^{n+1}-X_s^n\|^2_1+\|Y_s^{n+1}-Y_s^{n}\|^2)ds\right]\nonumber\\
&&+C\mathbb{E}\left[\int_0^t(\|X_s^{n}-X_s^{n-1}\|^2_1+\|Y_s^{n}-Y_s^{n-1}\|^2)ds\right]\label{Iterating}.
\end{eqnarray}
Iterating inequality \eqref{Iterating}  we obtain
\begin{eqnarray*}
\mathbb{E}\sup\limits_{0\leq s\leq
T_0}\left[\|X_s^{n+1}-X_s^{n}\|_1^2+\|Y_s^{n+1}-Y_s^{n}\|^2\right]\leq
C\frac{(C_{T_0})^n}{n!}.
\end{eqnarray*}
 This implies that there exists  an $\mathscr{F}_t-$adapted  couple
 $$(X_t, Y_t)\in L^2(\Omega; {C}([0, T_0]; \mathcal {H}^1))\times L^2(\Omega;
{C}([0, T_0]; H)),$$ which satisfies
\begin{eqnarray*}
\lim\limits_{n\rightarrow\infty}\mathbb{E}\sup\limits_{0\leq t\leq
T_0}\left[\|X_t^{n}-X_t\|_1^2+\|Y_t^{n}-Y_t\|^2\right]=0.
\end{eqnarray*}
Thanks to inequality \eqref{Energy-Slow-1-1}, we can show that
$\{\dot{X}^{n}_t\}_{n\in \mathbb{N}}$ is also a  Cauchy sequence
converging in  $L^2(\Omega; {C}([0, T_0]; H))$  to the limit
$\dot{X}_t$. Now, taking the limit in Eq.
\eqref{Nonlinear-Appro-Slow} and Eq. \eqref{Nonlinear-Appro-Fast},
it is seen that $(X_t, Y_t)$ is a solution to Eqs.
\eqref{Nonlinear-Integral-Slow}-\eqref{Nonlinear-Integral-Fast}. The
uniqueness is a directive consequence of energy equalities and
Gronwall's inequality.

Now, it remains to prove the energy equalities.  For all $ t\in [0,
T_0]$ one has the following convergence in mean-square as
$n\rightarrow\infty$:
\begin{eqnarray}
&&\dot{X}_t^{n}\rightarrow \dot{X}_t, \label{Conver-1}\\
&&X_t^n\rightarrow X_t.\label{Conver-2}
\end{eqnarray}
 {Also, we have the convergence }
\begin{eqnarray}
&&\int_0^t\big(X_s^n,
\sigma(X_s^n)dW^1_s\big)_Hds\rightarrow\int_0^t\big(X_s,
\sigma(X_s)dW^1_s\big)_Hds,\label{Conver-3}
\end{eqnarray}
\begin{eqnarray}
\int_0^t\|\sigma(X_s^n)\|^2_{Q_1}ds\rightarrow\int_0^t\|\sigma(X_s)\|^2_{Q_1}ds,
\end{eqnarray}
and
\begin{eqnarray}
\int_0^t\left(\dot{X}_s^n, f(X_s^n,
Y_s^n)\right)_Hds\rightarrow\int_0^t\left(\dot{X}_s, f(X_s,
Y_s)\right)_Hds\label{Conver-4}
\end{eqnarray}
in mean for all $0\leq t\leq T_0$, as $n$ goes to infinity.  Then,
by selecting a subsequence if necessary, we can assume  the
$\mathbb{P}-a.s.$ convergence in \eqref{Conver-1}-\eqref{Conver-4}.
As a result, one can obtain the energy equality given by
\eqref{energy-nonlinear-slow}. By a similar arguments we can get the
energy equality \eqref{energy-nonlinear-fast}.
\end{proof}
\end{theorem}

\section{Ergodicty for frozen equation}\label{section-ergodicty}
This section focuses on the ergodicty  for the fast motion equation.
For fixed $x\in H$ consider the problem associated with the fast
motion equation with a frozen slow component:
\begin{eqnarray}
&&\frac{\partial Y_t(\xi)}{\partial t}=\Delta Y_t(\xi)+g(x, Y_t(\xi))+{b}(x,Y_t(\xi))\dot{W_t^2},\label{Frozen-Fast}\\
&&Y_t(\xi)=0, (\xi, t)\in \partial D\times [0, \infty),\label{Frozen-Boundary}\\
&&Y_0(\xi)=y.\label{Frozen-Initial}
\end{eqnarray}
As the coefficients $g$ and $b$  are Lipschitz-continuous, for any
fixed slow component $x\in H$ and any initial value $y\in H$, the
system \eqref{Frozen-Fast}-\eqref{Frozen-Initial} admits a unique
mild solution denoted by $Y_t^{x, y}$. Moreover, we have the energy
equality
\begin{eqnarray*}
\|Y_t^{x, y}\|^2&=&\|y\|^2+2\int_0^t\langle AY_s^{x, y}, Y_s^{x,
y}\rangle ds+2\int_0^t\big(g(x,Y_s^{x, y}), Y_s^{x,
y}\big)_Hds\nonumber\\
&&+2\int_0^t\big(Y_s^{x, y},
b(x,Y_s^{x, y})dW_s^{2}\big)_H+\int_0^t\|b(x,Y_s^{x, y})\|^2_{Q_2}ds,\;\;a.s.,
\end{eqnarray*}
which implies the differential form
\begin{eqnarray}
\frac{d}{dt}\mathbb{E}\|Y_t^{x, y}\|^2&=&2\mathbb{E}\langle AY_t^{x,
y}, Y_t^{x, y}\rangle+2\mathbb{E}\big(g(x,Y_t^{x, y}), Y_t^{x,
y}\big)_H\nonumber\\
&&+\mathbb{E}\|b(x,Y_t^{x, y})\|^2_{Q_2},\label{Frozen-energy}
\end{eqnarray}
and then, thanks to the Poincar\'{e} inequality and the Lipschitz
continuity of $g$ and $b$, we obtain
\begin{eqnarray}
\frac{d}{dt}\mathbb{E}\|Y_t^{x,
y}\|^2&\leq&-2\alpha_1\mathbb{E}\|Y_t^{x, y}\|^2
+2\mathbb{E}\left|\big(g(x,Y_t^{x, y})-g(x,0), Y_t^{x,
y}\big)_H\right|\nonumber\\
&&+2\mathbb{E}\left|\big(g(x,0), Y_t^{x,
y}\big)_H\right|+\|b(x, Y_t^{x, y})-b(x, 0)\|^2_{Q_2}\nonumber\\
&&+ \|b(x, 0)\|^2_{Q_2}+2\Big|\Big(b(x, Y_t^{x, y})-b(x, 0),b(x, 0)\Big)_{Q_2}\Big|\nonumber\\
&\leq&-(2\alpha_1-2L_g- L^2_{b}-\rho)\mathbb{E}\|Y_t^{x,
y}\|^2+C_\rho(1+\|x\|^2)\nonumber
\end{eqnarray}
with positive constants $\rho$ and $C_\rho$ independent of $x$ and
$Y_t^{x,y}$, {where we have used the Young inequality in the form
$|a_1a_2|\leq \rho|a_1|^2+C_{\rho}|a_2|^2$ for $\rho>0$ at the
second step. Hence taking  \eqref{decay} into account and choosing
$\rho$ small enough we can find $C_1, C_2>0$ such that
\begin{eqnarray*}
\frac{d}{dt}\mathbb{E}\|Y_t^{x, y}\|^2&\leq&-C_1\mathbb{E}\|Y_t^{x,
y}\|^2+C_2(1+\|x\|^2),
\end{eqnarray*}
which, by Gronwall's inequality, allows us to get
\begin{equation}\label{Fast-motion-energy-bound}
\mathbb{E}\|Y_t^{x, y}\|^2\leq
C\left(e^{-ct}\|y\|^2+\|x\|^2+1\right),\; t>0
\end{equation}
for some constants $c,C>0$.

Now, for any $x\in H$,  let $\{P_t^x\}_{t\geq 0}$ denote the Markov
transition semigroup corresponding to the system
\eqref{Frozen-Fast}-\eqref{Frozen-Initial}, which is defined by
\begin{eqnarray*}
P_t^x\psi(y)=\mathbb{E}\psi(Y_t^{x, y}),\;t\geq0, y\in H
\end{eqnarray*}
for $\psi\in \mathcal {B}_b(H)$, the space of bounded measurable
functions on $H$. It follows from \eqref{Fast-motion-energy-bound}
that  there exists an invariant probability  measure $\mu^x$ for
$P^x_t$   on  $H$ such that
$$ \int_HP_t^x\psi(y) \mu^x(dy)=\int_H\psi(y)\mu^x(dy), \quad
t\geq 0
$$
for any  $\psi \in \mathcal {B}_b(H)$  (for a proof, see, e.g.,
Cerrai \cite{Cerrai2}, Section 2.1). Then by repeating the standard
argument as in the proof of Proposition 4.2 in Cerrai
\cite{Cerrai-Siam}, one can easily establish   the following
estimate:
\begin{equation}\label{mu-Momenent-bound}
\int_H\|y\|^2\mu^x(dy)\leq C(1+\|x\|^2).
\end{equation}
Let $Y^{x, y'}_t$ be the mild solution to
system\eqref{Frozen-Fast}-\eqref{Frozen-Initial} with the initial
value $Y_0=y'$. We can show that for any $t\geq0$,
\begin{eqnarray*}
\mathbb{E}\|Y_t^{x, y}-Y_t^{x, y'}\|^2\leq C\|y-y'\|^2e^{-c t}
\end{eqnarray*}
with $C,c>0,$ which allows us  to conclude that $\mu^x$ is the
unique invariant measure for $P^x_t$. Then, using the invariant
property of $\mu^x$, \eqref{mu-Momenent-bound} and condition (A1),
we obtain the following  inequality
\begin{eqnarray}
\nonumber\left\|\mathbb{E}f(x, Y_t^{x, y})-\int_Hf(x,
z)\mu^x(dz)\right\|^2&=&\left\|\int_H\mathbb{E}\big(f(x, Y_t^{x,
y})-f(x, Y_t^{x, z})\big)\mu^x(dz)\right\|^2\\
\nonumber&\leq&C\int_H\mathbb{E}\left\|Y_t^{x, y}-Y_t^{x,
z}\right\|^2\mu^x(dz)\\
\nonumber&\leq&Ce^{-c t}\int_H\|y-z\|^2\mu^x(dz)\\
&\leq&Ce^{-c
t}\big(1+\|x\|^2+\|y\|^2\big).\label{Averaging-Expectation}
\end{eqnarray}

\section{A priori bounds and auxiliary  processes}
We now prove the following estimates for the solution to the system
\eqref{Integral-Slow}-\eqref{Integral-Fast}.
\begin{lemma}\label{5-1}
Assume that $X_0\in \mathcal {H}^1, \dot{X}_0\in H, Y_0\in H$ and
(A1)-(A5) hold, then there exists a constant $C_{T_0}>0$ such that
\begin{equation}\label{Boundness-slow-solution}
\sup\limits_{\begin{subarray}{l}
         \;\:\,\epsilon>0, \\
         0\leq t\leq T_0\\
    \end{subarray}} \mathbb{E}(\|X^\epsilon_t\|_1^2+\|\dot{X}^\epsilon_t\|^2)\leq
C_{T_0}\left(1+\|Y_0\|^2+\|X_0\|^2_1+\|\dot{X}_0\|^2\right)
\end{equation}
and
\begin{equation}\label{Boundness-fast-solution}
\sup\limits_{\begin{subarray}{l}
        \epsilon>0, \\
         0\leq t\leq T_0\\
    \end{subarray}} \mathbb{E}\|Y^\epsilon_t\|^2\leq
    C_{T_0}\left(1+\|Y_0\|^2+\|X_0\|^2_1+\|\dot{X}_0\|^2\right).
\end{equation}
\begin{proof}
By the energy equality \eqref{energy-nonlinear-slow} and sublinear
growth condition of $f$, we have
\begin{eqnarray}
\mathbb{E}\left(\|\dot{X}_t^\epsilon\|^2+\|X_t^\epsilon\|^2_1\right)&=&\|\dot{X}_0\|^2+\|X_0\|_1^2+2\mathbb{E}\int_0^t\big(\dot{X}_s^\epsilon,
f(X_s^\epsilon, Y_s^\epsilon)\big)_Hds \nonumber\\
&&+\int_0^t\mathbb{E}\|\sigma(X_s^\epsilon)\|^2_{Q_1}ds\nonumber\\
&\leq&\|\dot{X}_0\|^2+\|X_0\|_1^2+C\int_0^t\mathbb{E}\left(\|\dot{X}_s^\epsilon\|^2+\|X_s^\epsilon\|^2_1\right)ds\nonumber\\
&&+C\int_0^t\mathbb{E}\left(1+\|Y_s^\epsilon\|^2\right)ds,\nonumber
\end{eqnarray}
so that
\begin{eqnarray}
\mathbb{E}\left(\|\dot{X}_t^\epsilon\|^2+\|X_t^\epsilon\|^2_1\right)&\leq&
e^{Ct}\left(\|\dot{X}_0\|^2+\|X_0\|_1^2\right)\nonumber\\
&&+C\int_0^te^{C(t-s)}\left(1+\mathbb{E}\|Y_s^\epsilon\|^2\right)ds.\label{slow
-moment-bounds}
\end{eqnarray}
 Thanks to energy equality \eqref{energy-nonlinear-fast} {and  condition (A5)}, we
 have the ordinary differential inequality
\begin{eqnarray}
\frac{d}{dt}\mathbb{E}\|Y^\epsilon_t\|&=&\frac{2}{\epsilon}\mathbb{E}\langle
AY^\epsilon_t, Y^\epsilon_t\rangle
+\frac{2}{\epsilon}\mathbb{E}\big(g(X^\epsilon_t, Y^\epsilon_t),
Y^\epsilon_t\big)_H +\frac{1}{\epsilon}\mathbb{E}\|b(X_t^\epsilon,Y_t^\epsilon)\|^2_{Q_2}ds\nonumber\\
&\leq&-\frac{C_1}{\epsilon}\mathbb{E}\|Y^\epsilon_t\|^2
+\frac{C_2}{\epsilon}\left(1+\mathbb{E}\|X_t^\epsilon\|^2\right),\nonumber
\end{eqnarray}
and hence
\begin{eqnarray*}
\mathbb{E}\|Y^\epsilon_t\|\leq
e^{C_2t}\|Y_0\|^2+\frac{C_2}{\epsilon}\int_0^te^{-\frac{C_1}{\epsilon}(t-s)}\left(1+\mathbb{E}\|X_s^\epsilon\|^2\right)ds.
\end{eqnarray*}
According to \eqref{slow -moment-bounds}, we obtain
\begin{eqnarray*}
\mathbb{E}\|Y^\epsilon_t\|&\leq&
C_{T_0}\left(1+\|Y_0\|^2+\|\dot{X}_0\|^2+\|X_0\|_1^2\right)\nonumber\\
&&+\frac{C_{T_0}}{\epsilon}\int_0^te^{-\frac{C_1}{\epsilon}(t-s)}\int_0^s\mathbb{E}\|Y_r^\epsilon\|^2drds,\;t\in
[0, T_0].
\end{eqnarray*}
By change of variables, this yields
\begin{eqnarray*}
\mathbb{E}\|Y^\epsilon_t\|&\leq&
C_{T_0}\left(1+\|Y_0\|^2+\|\dot{X}_0\|^2+\|X_0\|_1^2\right)\nonumber\\
&&+C_{T_0}\int_0^t\mathbb{E}\|Y_r^\epsilon\|^2\left[\int_0^{\frac{t-r}{\epsilon}}e^{-C_1v}dv\right]dr\nonumber\\
&\leq&C_{T_0}\left(1+\|Y_0\|^2+\|\dot{X}_0\|^2+\|X_0\|_1^2\right)\nonumber\\
&&+C_{T_0}\int_0^t\mathbb{E}\|Y_r^\epsilon\|^2dr.
\end{eqnarray*}
Hence, by Gronwall's inequality we get
\begin{eqnarray*}
\mathbb{E}\|Y^\epsilon_t\|&\leq&
C_{T_0}\left(1+\|Y_0\|^2+\|\dot{X}_0\|^2+\|X_0\|_1^2\right),
\end{eqnarray*}
which, gives the estimate \eqref{Boundness-fast-solution}. By
replacing the above estimate  in \eqref{slow -moment-bounds} and
using the Gronwall inequality again, we obtain the   estimate
\eqref{Boundness-slow-solution}.
\end{proof}
\end{lemma}
We now provide a regularity estimate of $X^\epsilon_t$ in the time
variable.

\begin{lemma}\label{5-2}
Suppose that conditions in Lemma \ref{5-1} hold. For any $X_0\in
\mathcal {H}^1, \dot{X}_0\in H, Y_0\in H$, {there exists a constant}
$C_{T_0}>0$ such that for any $0\leq t\leq t+h\leq T_0$,
\begin{equation}\label{Regularity}
\sup\limits_{\epsilon>0
}\mathbb{E}\|X^\epsilon_{t+h}-X^\epsilon_{t}\|^2\leq
C_{T_0}(\|X_0\|_1^2+\|\dot{X}_0\|^2+\|Y_0\|^2)h^2.
\end{equation}
\begin{proof}
Clearly, we have
\begin{eqnarray*}
\mathbb{E}\|X^{\epsilon}_{t+h}-X^\epsilon_{t}\|^2&=&\mathbb{E}\|\int_t^{t+h}\dot{X}^\epsilon_{s}ds\|^2\\
&\leq&h\int_t^{t+h}\mathbb{E}\|\dot{X}^\epsilon_{s}\|^2ds,
\end{eqnarray*}
so that, by  referring to Lemma\eqref{Boundness-slow-solution}, we
have
\begin{eqnarray*}
\mathbb{E}\|X^{\epsilon}_{t+h}-X^\epsilon_{t}\|^2\leq
C_{T_0}(\|X_0\|_1^2+\|\dot{X}_0\|^2+\|Y_0\|^2)h^2.
\end{eqnarray*}
This completes the proof.
\end{proof}
\end{lemma}

Next, we define two auxiliary  processes. To this end, divide the
time interval $[0, T_0]$ into  subintervals of equal length
$\delta$. We construct a process $\hat{Y}_t^\epsilon$, with initial
value $\hat{Y}_0^\epsilon=Y_0$, which is described by equation
\begin{eqnarray*}
d\hat{Y}_t^\epsilon=\frac{1}{\epsilon}A\hat{Y}_t^\epsilon
dt+\frac{1}{\epsilon}g(X_{k\delta}^\epsilon,
\hat{Y}_t^\epsilon)dt+\frac{1}{\sqrt{\epsilon}}b(X_{k\delta}^\epsilon,
\hat{Y}_t^\epsilon)dW_t^{2}
\end{eqnarray*}
for $t\in \left[k\delta, \min\big\{(k+1)\delta, T_0\big\}\right ),
k\geq0$, where $X_{k\delta}^\epsilon$ is the slow solution process
at time $k\delta$. For the left end of each subinterval we set
$$\hat{Y}_{(k+1)\delta}^\epsilon=\lim\limits_{t\rightarrow(k+1)\delta-}\hat{Y}_t^\epsilon.$$
Denote $\lfloor\cdot\rfloor$ to be the integer function.  The mild
form for $\hat{Y}_t^\epsilon$ takes form
\begin{eqnarray}
&&\hat{Y}^\epsilon_t=G_{t/\epsilon}Y_0+\frac{1}{\epsilon}\int_0^tG_{(t-s)/\epsilon}g(X^\epsilon_{s(\delta)},
\hat{Y}_s^\epsilon)ds\nonumber\\
&&\quad\quad\;+\frac{1}{\sqrt{\epsilon}}\int_0^tG_{(t-s)/\epsilon}b(X^\epsilon_{s(\delta)},
\hat{Y}_s^\epsilon)dW^2_s,\;t\in [0, T_0],\nonumber
\end{eqnarray}
where $s(\delta)=\lfloor s/\delta\rfloor\delta$ is the nearest
breakpoint preceding $s$. Define the process $\hat{X}_t^\epsilon$ by
integral
\begin{eqnarray*}
\hat{X}_t^\epsilon=S_t'X_0+S_t\dot{X}_0+\int_0^tS_{t-s}f(X^\epsilon_{s(\delta)},
\hat{Y}_s^\epsilon)ds+\int_0^tS_{t-s}\sigma(X_{s(\delta)}^\epsilon)dW^1_s.
\end{eqnarray*}
According to the previous lemma, we can easily prove the following
boundedness of the second moment of $\hat{Y}_t^\epsilon$.
\begin{lemma}\label{Y-hat}
Suppose that conditions in Lemma \ref{5-1} hold.  Then there exists
a constant $C_{T_0}>0$ such that for any $t\in [0, T_0]$ it holds
\begin{equation}\label{Y-hat-1}
\mathbb{E}\|\hat{Y}_t^\epsilon\|^2\leq C_{T_0}.
\end{equation}
\begin{proof}
For $t\in [0, T_0]$ 
we have
\begin{eqnarray*}
\|\hat{Y}_t^\epsilon\|^2&=&\|{Y}_0\|^2+\frac{2}{\epsilon}\int_{0}^t\langle
A\hat{Y}_s^\epsilon, \hat{Y}_s^\epsilon\rangle
ds+\frac{2}{\epsilon}\int_{0}^t\big(g(X_{s(\delta)}^\epsilon,
\hat{Y}_s^\epsilon),\hat{Y}_s^\epsilon\big)_Hds\\
&+&\frac{2}{\epsilon}\int_{0}^t\big(\hat{Y}_s^\epsilon,dW_s^2\big)_H+\frac{1}{\epsilon}\int_{0}^t\|b(X_{s(\delta)}^\epsilon,
\hat{Y}_s^\epsilon)\|^2_{Q_2}ds,\;a.s..
\end{eqnarray*}
This implies a differential equality in the form
\begin{eqnarray*}
\frac{d}{dt}\mathbb{E}\|\hat{Y}_t^\epsilon\|^2=\frac{2}{\epsilon}\mathbb{E}\langle
A\hat{Y}_t^\epsilon,
\hat{Y}_t^\epsilon\rangle+\frac{2}{\epsilon}\mathbb{E}\big(g(X_{t(\delta)}^\epsilon,
\hat{Y}_t^\epsilon),\hat{Y}_t^\epsilon\big)_H+\frac{1}{\epsilon}\mathbb{E}\|b(X_{t(\delta)}^\epsilon,
\hat{Y}_t^\epsilon)\|^2_{Q_2}.
\end{eqnarray*}
Using similar {arguments for \eqref{Frozen-energy}}, we obtain
\begin{eqnarray*}
\frac{d}{dt}\mathbb{E}\|\hat{Y}_t^\epsilon\|^2\leq-\frac{1}{\epsilon}C_1\mathbb{E}\|\hat{Y}_t^\epsilon\|^2
+\frac{1}{\epsilon}C_2(1+\mathbb{E}\|X^\epsilon_{t(\delta)}\|^2).
\end{eqnarray*}
Then, due to \eqref{Boundness-slow-solution} we have
\begin{eqnarray*}
\frac{d}{dt}\mathbb{E}\|\hat{Y}_t^\epsilon\|^2\leq-\frac{1}{\epsilon}C_1\mathbb{E}\|\hat{Y}_t^\epsilon\|^2
+\frac{1}{\epsilon}C_2
\end{eqnarray*}
and the proof is completed by an application of Gronwall's
inequality.
\end{proof}
\end{lemma}

The following difference estimates will be used in the proof of our
strong error result stated in Theorem \ref{Theorem} below.
\begin{lemma}\label{main-lemma}
Suppose that conditions in Lemma \ref{5-1} hold. Then there exists a
constant $C>0$ such that for any $t\in [0, T_0]$ it holds
\begin{equation}\label{Auxiliary-fast-difference}
\mathbb{E}\|Y_t^\epsilon-\hat{Y}_t^\epsilon\|^2\leq C\delta^2,
\end{equation}
\begin{equation}\label{Auxiliary-slow-difference}
\mathbb{E}\|X_t^\epsilon-\hat{X}_t^\epsilon\|_1^2\leq C\delta^2
\end{equation}
and
\begin{equation}\label{Auxiliary-slow-derive}
\mathbb{E}\|\dot{X}_t^\epsilon-\dot{\hat{X}}_t^\epsilon\|^2\leq
C\delta^2.
\end{equation}
\begin{proof}
For $t\in [0, T_0]$ with $t\in [k\delta, (k+1)\delta)$ {we have}
\begin{eqnarray*}
\mathbb{E}\|Y_t^\epsilon-\hat{Y}_t^\epsilon\|^2&=&\mathbb{E}\|Y_{k\delta
}^\epsilon-\hat{Y}_{k\delta}^\epsilon\|^2+\frac{2}{\epsilon}\int_{k\delta}^t\mathbb{E}\langle
A(Y_s^\epsilon-\hat{Y}_s^\epsilon),
Y_s^\epsilon-\hat{Y}_s^\epsilon\rangle
ds\\
&&+\frac{2}{\epsilon}\int_{k\delta}^t\mathbb{E}\big(g(X_{k\delta}^\epsilon,
\hat{Y}_s^\epsilon)-g(X_s^
\epsilon,Y_s^\epsilon),Y_s^\epsilon-\hat{Y}_s^\epsilon\big)_Hds\\
&&+\frac{1}{\epsilon}\int_{k\delta}^t\mathbb{E}\|b(X^\epsilon_{k\delta},\hat{Y}^\epsilon_s)-b(X_s^
\epsilon,Y_s^\epsilon)\|^2_{Q_2}ds.
\end{eqnarray*}
This shows that
\begin{eqnarray}
\frac{d}{dt}\mathbb{E}\|Y_t^\epsilon-\hat{Y}_t^\epsilon\|^2&=&\frac{2}{\epsilon}
\mathbb{E}\langle A(Y_t^\epsilon-\hat{Y}_t^\epsilon),
Y_t^\epsilon-\hat{Y}_t^\epsilon\rangle\nonumber\\
&&+\frac{2}{\epsilon}
\mathbb{E}\big(g(X^\epsilon_{k\delta},\hat{Y}^\epsilon_t)-g(X_t^
\epsilon,Y_t^\epsilon),Y_t^\epsilon-\hat{Y}_t^\epsilon\big)_H\nonumber\\
&&+\frac{1}{\epsilon}
\mathbb{E}\|b(X^\epsilon_{k\delta},\hat{Y}^\epsilon_t)-b(X_t^
\epsilon,Y_t^\epsilon)\|^2_{Q_2},\label{4-5}
\end{eqnarray}
which, with the aid of the Young  inequality in the form
$|ab|\leq\frac{\rho}{2C}|a|^2+C_{\rho}|b|^2$  for $\rho>0$, yields
from (A5) that

\begin{eqnarray*}
\nonumber\frac{d}{dt}\mathbb{E}\|Y_t^\epsilon-\hat{Y}_t^\epsilon\|^2
&\leq&-
\frac{C_1}{\epsilon}\mathbb{E}\|Y_t^\epsilon-\hat{Y}_t^\epsilon\|^2
+\frac{C_2}{\epsilon}\mathbb{E}\|X_{k\delta}^\epsilon-X_t^\epsilon\|^2.
\end{eqnarray*}
By applying \eqref{Regularity}, we have
\begin{eqnarray*}
\frac{d}{dt}\mathbb{E}\|Y_t^\epsilon-\hat{Y}_t^\epsilon\|^2 &\leq&-
\frac{C_1}{\epsilon}\mathbb{E}\|Y_t^\epsilon-\hat{Y}_t^\epsilon\|^2+C_2\frac{\delta^2}{\epsilon}
\end{eqnarray*}
and then, due to {Gronwall's inequality}, we have
\begin{eqnarray}
\mathbb{E}\|Y_t^\epsilon-\hat{Y}_t^\epsilon\|^2\leq
e^{-C_1\frac{(t-k\delta)}{\epsilon}}\mathbb{E}\|Y_{k\delta
}^\epsilon-\hat{Y}_{k\delta}^\epsilon\|^2+C_2(1-e^{-C_1\frac{(t-k\delta)}{\epsilon}})\delta^2.\label{5-10}
\end{eqnarray}
Taking $t=(k+1)\delta$ in above inequality, one obtains
\begin{equation*}
\mathbb{E}\|Y_{(k+1)\delta}^\epsilon-\hat{Y}_{(k+1)\delta}^\epsilon\|^2\leq\mathbb{E}\|Y_{k\delta}^\epsilon-\hat{Y}_{k\delta}^\epsilon\|^2
e^{-C_1\frac{\delta}{\epsilon}}+C_2\left(1-e^{-C_1\frac{\delta}{\epsilon}}\right)\delta^2.
\end{equation*}
Iterating the above inequality recursively from $k$ to $0$, we get
\begin{eqnarray*}
\mathbb{E}\|Y_{(k+1)\delta}^\epsilon-\hat{Y}_{(k+1)\delta}^\epsilon\|^2
&\leq& C\left(1-e^{-C_1\frac{\delta}{\epsilon}}\right)\delta^2
\sum\limits_{0\leq j\leq k}e^{-jC_1\frac{\delta}{\epsilon}}\nonumber\\
&\leq&C\delta^2.
\end{eqnarray*}
Then, by using again \eqref{5-10}, we have
\begin{equation*}
\mathbb{E}\|Y_t^\epsilon-\hat{Y}_t^\epsilon\|^2\leq C\delta^2,
\end{equation*}
and therefore the first assertion follows.

For the second part,  we employ the inequality \eqref{S-t-2} in
{Proposition} \ref{lemma-2}  and \eqref{S-t-1-stocha} in Proposition
\ref{Pro-2.3} to obtain
 \begin{eqnarray*}
\mathbb{E}\|X_t^\epsilon-\hat{X}_t^\epsilon\|^2_1
&\leq&C\int_0^{T_0}\mathbb{E}\|X_s^\epsilon-{X}_{s(\delta)}^\epsilon\|^2ds
+C\int_0^{T_0}\mathbb{E}\|Y_s^\epsilon-\hat{Y}_s^\epsilon\|^2ds.
\end{eqnarray*}
By taking \eqref{Regularity} and the first estimate
\eqref{Auxiliary-fast-difference} into account, this yields the
desired  estimate \eqref{Auxiliary-slow-difference}.

For the remaining one,  by \eqref{S-t-3} in {Proposition}
\ref{lemma-2}  and \eqref{S-t-2-stocha} in Proposition
\ref{Pro-2.3}, we have
\begin{eqnarray}
\mathbb{E}\|\dot{X}_t^\epsilon-\dot{\hat{X}}_t^\epsilon\|^2&\leq&
T_0\int_0^{T_0}\mathbb{E}\|f(X^\epsilon_{s(\delta)},\hat{Y}^\epsilon_s)-f(X_s^\epsilon,Y^\epsilon_s)\|^2ds\nonumber\\
&&+\nonumber C\int_0^{T_0}\|\sigma(X^\epsilon_{s(\delta)})-\sigma(X_s^\epsilon)\|^2_{Q_1}ds\\
&\leq&
C\int_0^{T_0}\mathbb{E}\left(\|X^\epsilon_{s(\delta)}-X_s^\epsilon\|^2+\|\hat{Y}^\epsilon_s-Y^\epsilon_s\|^2\right)ds.\nonumber
\end{eqnarray}
The combination of  \eqref{Regularity} and
\eqref{Auxiliary-fast-difference} yields
\begin{eqnarray}
\mathbb{E}\|\dot{X}_t^\epsilon-\dot{\hat{X}}_t^\epsilon\|^2\leq
C\delta^2.\nonumber
\end{eqnarray}
This completes the proof.
\end{proof}
\end{lemma}

\section{Main result}
In this section we  will consider the  effective dynamics system
\begin{eqnarray}
&&\label{Averaing-equation}\frac{\partial^2\bar{X}_t(\xi)}{\partial
t^2}=\Delta \bar{X}_t(\xi)+ \bar{f}(\bar{X}_t(\xi))+\sigma(\bar{X}_t(\xi))\dot{W}_t^{1},\\
&&\bar{X}_t(\xi)=0, (\xi, t)\in \partial
D\times [0, \infty),\label{Averaing-equation-Boundary}\\
&&\bar{X}_0(\xi)=X_0(\xi), \frac{\partial \bar{X}_t(\xi)}{\partial
t}\Big|_{t=0}=\dot{X}_0(\xi), \;\xi\in
D,\label{Averaing-equation-inital-data}
\end{eqnarray}
with
$$
\bar{f}(x)=\int_Hf(x, y)\mu^x(dy),\quad x\in H,
$$
where $\mu^x$ denotes the unique invariant measure for system
\eqref{Frozen-Fast}-\eqref{Frozen-Initial} introduced in Section
\ref{section-ergodicty}. Moreover, due to Cerrai and Freidlin
\cite{Cerrai1}, the mapping $\bar{f}: H\mapsto H$ is Lipschitz
continuous. Furthermore, by taking \eqref{mu-Momenent-bound} into
account, we have
\begin{equation}
\|f(x, y)-\bar{f}(x)\|^2\leq
C_f\left(1+\|x\|^2+\|y\|^2\right)\label{Averaging-Function},\;\;x,y
\in H.
\end{equation}
The mild form for the averaged system
\eqref{Averaing-equation}-\eqref{Averaing-equation-inital-data} is
given by
\begin{eqnarray}
\bar{X}_t=S_t'X_0+S_t\dot{X}_0+\int_0^tS_{t-s}\bar{f}(\bar{X}_s
)ds+\int_0^tS_{t-s}\sigma(\bar{X}_s
)dW^1_s.\label{Integral-Averaging-equation}
\end{eqnarray}
Since we are working with a global Lipschitz nonlinearity, standard
arguments  guarantee the existence of a unique mild solution to the
above integral equation in $L^2(\Omega; {C}([0, T_0]; \mathcal
{H}^1))$ with its time derivative $\dot{\bar{X}}_t\in L^2(\Omega;
{C}([0, T_0]; H))$.

Now, we present averaging result  with explicit error bound, which
reveals that the order of strong convergence in averaging is
$\frac{1}{2}$.

\begin{theorem}\label{Theorem}
Let assumptions (A1)-(A5) be fulfilled.  Then there exists a
constant $C$ such that for any $t\in [0, T_0]$,
\begin{eqnarray*}
\mathbb{E}\left
(\|{X}_t^\epsilon-\bar{X}_t\|^{2}_1+\|\dot{{X}}_t^\epsilon-\dot{\bar{X}}_t\|^{2}\right
)\leq C{\epsilon},
\end{eqnarray*}
here $\bar{X}_t$  is the mild solution of the effective system
\eqref{Averaing-equation}-\eqref{Averaing-equation-inital-data}.
\end{theorem}
The proof of Theorem \ref{Theorem} is postponed to the end of this
section.  {First, we   give some lemmas that will be used.}
\begin{lemma}\label{Key-lemma}
Suppose that conditions in Theorem \ref{Theorem} hold. Then there
exists a constant $C$ such that for any $t\in [0, T_0]$,
\begin{equation*}
\mathbb{E}\left
(\|\hat{X}_t^\epsilon-\bar{X}_t\|^{2}_1+\|\dot{\hat{X}}_t^\epsilon-\dot{\bar{X}}_t\|^{2}\right
)\leq C({\epsilon}+\delta^2).
\end{equation*}
\begin{proof}
For any $t\in [0, T_0]$ we have the decomposition
\begin{equation}
\hat{X}_t^\epsilon-\bar{X}_t=\mathcal {L}_t^\epsilon+\mathcal
{M}_t^\epsilon+\mathcal{N}_t^\epsilon+\mathcal
{Q}_t^\epsilon+\mathcal {T}_t^\epsilon,\label{Devided}
\end{equation}
where
\begin{eqnarray*}
&&\mathcal
{L}_t^\epsilon=\int_0^tS_{t-s}\big(f(X^\epsilon_{s(\delta)},
\hat{Y}_s^\epsilon)-\bar{f}(X_s^\epsilon)\big)ds,\label{A-1}\\
&&\mathcal
{M}_t^\epsilon=\int_0^tS_{t-s}\big(\bar{f}(X_s^\epsilon)-\bar{f}(\hat{X}_s^\epsilon)\big)ds,\label{A-2}  \\
&&\mathcal{N}_t^\epsilon=\int_0^tS_{t-s}\big(\bar{f}(\hat{X}_s^\epsilon)-\bar{f}(\bar{X}_s)\big)ds,\label{A-3} \\
&&\mathcal
{Q}_t^\epsilon=\int_0^tS_{t-s}\big(\sigma(X_s^\epsilon)-\sigma(\hat{X}_s^\epsilon)\big)dW_s^1\label{A-4}
\end{eqnarray*}
and
\begin{eqnarray*}
\mathcal
{T}_t^\epsilon=\int_0^tS_{t-s}\big(\sigma(\hat{X}_s^\epsilon)-\sigma(\bar{X}_s^\epsilon)\big)dW_s^1.\label{A-5}
\end{eqnarray*}
By virtue of  Proposition \ref{lemma-2}, Proposition \ref{Pro-2.3}
and \eqref{Auxiliary-slow-difference}, it is easy to show that for
any $t\in[0, T_0],$
\begin{eqnarray}
\mathbb{E}\|\mathcal {M}_t^\epsilon\|^2_1+\mathbb{E}\|\dot{\mathcal
{M}}_t^\epsilon\|^2+\mathbb{E}\|\mathcal
{Q}_t^\epsilon\|^2_1+\mathbb{E}\|\dot{\mathcal
{Q}}_t^\epsilon\|^2&\leq&
T_0\int_0^{T_0}\mathbb{E}\|X_s^\epsilon-\hat{X}_s^\epsilon\|^2ds\nonumber\\
&\leq&C\int_0^{T_0}\mathbb{E}\|X_s^\epsilon-\hat{X}_s^\epsilon\|^2_1ds\nonumber \\
&\leq& C\delta^2, \label{M-t-1}
\end{eqnarray}
where we used the fact that $\|\cdot\|$-norm can be dominated by
$\|\cdot\|_1$-norm at the second step. By using again the
Proposition \ref{lemma-2} and Proposition \ref{Pro-2.3}, we have
\begin{eqnarray}
\mathbb{E} \|\mathcal {N}_t^\epsilon\|^2_1+\mathbb{E}
\|\dot{\mathcal {N}}_t^\epsilon\|^2+\mathbb{E}\|\mathcal
{T}_t^\epsilon\|^2_1+\mathbb{E}\|\dot{\mathcal {T}}_t^\epsilon\|^2
&\leq&
T_0\int_0^{t}\mathbb{E}\|\hat{X}_s^\epsilon-\bar{X}_s^\epsilon\|^2ds\nonumber\\
&\leq& C\int_0^{t}\mathbb{E}
\|\hat{X}_s^\epsilon-\bar{X}_s^\epsilon\|^2_1ds,\label{N-t-1}
\end{eqnarray}
Now, due to the inequalities \eqref{M-t-1}, \eqref{N-t-1} and
\eqref{sup-A-1} in Lemma \ref{6-2} below, we can get
\begin{eqnarray*}
&&\mathbb{E} \left(\|\hat{X}_t^\epsilon-\bar{X}_t^\epsilon\|^2_1+\|\dot{\hat{X}}_t^\epsilon-\dot{\bar{X}}_t^\epsilon\|^2\right)\\
&&\leq C(
 {\epsilon+\delta^2} ) +C\int_0^{t}\mathbb{E}
 \|\hat{X}_s^\epsilon-\bar{X}_s^\epsilon\|^2_1ds,
\end{eqnarray*}
and the proof is complete in view of the Gronwall inequality.
\end{proof}
\end{lemma}
Now, we prove the following estimate, which is crucial in analyzing
the rate of strong error for averaging approximation.
\begin{lemma}\label{6-2}
Suppose that conditions in Lemma \ref{Key-lemma} hold. Then for any
$t\in [0,T_0 ] $, we have
\begin{eqnarray}\label{sup-A-1}
\mathbb{E}\left(\|\mathcal{L}_t^\epsilon\|^2_1+\|\dot{\mathcal{L}_t^\epsilon}\|^2\right)\leq
C(\delta^2+ {\epsilon} ),
\end{eqnarray}
where $C$ is a constant  independent of $(\epsilon, \delta).$

\begin{proof}
The proof will be divided into several steps. For any $t\in [0,
T_0)$, there exists an $n_t=\lfloor t/\delta\rfloor$ such that
$t\in[n_t\delta, (n_t+1)\delta\wedge T_0)$. Therefore, we have the
representation  in the form
\begin{equation}\label{I_1--I_3}
\mathcal {L}^\epsilon_t= I_1({t}, \epsilon)+I_2({t},
\epsilon)+I_3({t}, \epsilon),
\end{equation}
where
\begin{eqnarray*}
&&I_1({t}, \epsilon)=\sum\limits_{k=0}^{\lfloor
t/\delta\rfloor-1}\int_{k\delta}^{(k+1)\delta}
S_{t-s}\left(f(X_{k\delta}^\epsilon,
\hat{Y}_s^\epsilon)-\bar{f}(X_{k\delta}^\epsilon)\right)ds,\\
&&I_2({t}, \epsilon)=\sum\limits_{k=0}^{\lfloor
t/\delta\rfloor-1}\int_{k\delta}^{(k+1)\delta}
S_{t-s}\left(\bar{f}(X_{k\delta}^\epsilon)
-\bar{f}(X_s^\epsilon)\right)ds\\
&&\qquad\quad=\int_0^{\lfloor
t/\delta\rfloor\delta}S_{t-s} \left(\bar{f}(X_{s(\delta)}^\epsilon)-\bar{f}(X_s^\epsilon)\right)ds,\\
&&I_3({t}, \epsilon)=\int_{{\lfloor t/\delta\rfloor}\delta}^{t}
S_{t-s}\left(f(X_{{\lfloor t/\delta\rfloor}\delta}^\epsilon,
\hat{Y}^\epsilon_s)-\bar{f}(X_s^\epsilon)\right)ds.
\end{eqnarray*}

\textbf{Step 1:} Let us first deal with $I_2(t, \epsilon)$. Due to
the Lipschitz continuity of $\bar{f}$, we have the inequalities
\begin{eqnarray}
\|I_2(t, \epsilon)\|^2_1&\leq& \left[\int_0^{{\lfloor
t/\delta\rfloor}\delta}\left\|S_{t-s}\left(\bar{f}(X_{s(\delta)}^\epsilon)-\bar{f}(X_s^\epsilon)\right)\right\|_1ds\right]^2\nonumber\\
&\leq&T_0\int_0^{T_0}\left\|\bar{f}(X_{s(\delta)}^\epsilon)-\bar{f}(X_s^\epsilon)\right\|^2ds\nonumber\\
&\leq&C\int_0^{T_0}\left\|X_{s(\delta)}^\epsilon-X_s^\epsilon\right\|^2ds,\nonumber
\end{eqnarray}
so that, by Lemma \ref{5-2},
\begin{eqnarray*}
\mathbb{E} \|I_2(t, \epsilon)\|^2_1\leq
C_{T_0}\delta^2.
\end{eqnarray*}

\textbf{Step 2:} We proceed  to the estimate for $I_3(t, \epsilon)$.
As the mappings $f : H\times H\rightarrow H$ and $\bar{f}:
H\rightarrow H$ satisfy  sublinear growth condition, due to the
H\"{o}lder inequality we obtain
\begin{eqnarray}
\|I_3({t}, \epsilon)\|^2_1&\leq& \delta \int_{{\lfloor
t/\delta\rfloor}\delta}^{t} \left\|S_{t-s}\left(f(X_{{\lfloor
t/\delta\rfloor}\delta}^\epsilon,
\hat{Y}^\epsilon_s)-\bar{f}(X_s^\epsilon)\right)\right\|^2_1ds \nonumber\\
&\leq&\delta C\int_{\lfloor
t/\delta\rfloor\delta}^{t}\left(1+\|X_s^\epsilon\|^2+\|X_{{\lfloor
t/\delta\rfloor}\delta}^\epsilon\|^2+\|\hat{Y}^\epsilon_s\|^2\right)ds.\nonumber
\end{eqnarray}
With the aid of \eqref{Boundness-slow-solution} and \eqref{Y-hat-1},
we get
\begin{eqnarray*}
\mathbb{E}\|I_3({t}, \epsilon)\|^2_1&\leq& C\delta^2.
\end{eqnarray*}

\textbf{Step 3:}  { Let us now deal }with $I_1(t, \epsilon)$.
Concerning $I_1(t, \epsilon)$, we have by the series representation
of the Green's function that
\begin{eqnarray*}
I_1(t, \epsilon)=\sum\limits_{k=0}^{\lfloor
t/\delta\rfloor-1}\mathbbm{i}(t,k,\epsilon),
\end{eqnarray*}
where

\begin{eqnarray*}
&&\!\!\!\!\!\!\!\!\!\!\mathbbm{i}(t,k,\epsilon)\\
&=&\int_{k\delta}^{(k+1)\delta}\sum\limits_{i=1}^
\infty\frac{\sin\{\sqrt{\alpha_i}(t-s)\}}{\sqrt{\alpha_i}}\left(f(X_{k\delta}^\epsilon,
\hat{Y}_s^\epsilon)-\bar{f}(X_{k\delta}^\epsilon), e_i\right)_H\cdot
e_i ds\\
&=&\sum\limits_{i=1}^
\infty\frac{\sin\{{\sqrt{\alpha_i}t}\}}{\sqrt{\alpha_i}}\cdot
e_i\cdot \int_{k\delta}^{(k+1)\delta}\cos\{\sqrt{\alpha_i}
s\}\left(f(X_{k\delta}^\epsilon,
\hat{Y}_s^\epsilon)-\bar{f}(X_{k\delta}^\epsilon), e_i\right)_H
ds\\
&&-\sum\limits_{i=1}^
\infty\frac{\cos{\{\sqrt{\alpha_i}t\}}}{\sqrt{\alpha_i}}\cdot
e_i\cdot \int_{k\delta}^{(k+1)\delta}\sin\{\sqrt{\alpha_i}
s\}\left(f(X_{k\delta}^\epsilon,
\hat{Y}_s^\epsilon)-\bar{f}(X_{k\delta}^\epsilon), e_i\right)_H ds
\end{eqnarray*}
for $k=0,1,\cdots, \lfloor t/\delta\rfloor$.  Clearly, we have
\begin{eqnarray}
\|I_1(t, \epsilon)\|^2_1&=&\sum\limits_{k=0}^{\lfloor
t/\delta\rfloor-1}\|\mathbbm{i}(t,k,\epsilon)\|^2_1+2\sum\limits_{0\leq
i< j\leq \lfloor t/\delta\rfloor-1}\Big\langle
\mathbbm{i}(t,i,\epsilon),\mathbbm{i}(t,j,\epsilon)\Big\rangle_1\nonumber\\
&:=&\mathcal{A}_1(t,\epsilon)+2\mathcal{A}_2(t,\epsilon).\label{I_1-inner}
\end{eqnarray}
Note that
\begin{eqnarray*}
\|\mathbbm{i}(t,k,\epsilon)\|_1^2&\leq&2\sum\limits_{i=1}^
\infty\left[\int_{k\delta}^{(k+1)\delta}\cos\{\sqrt{\alpha_i}
s\}\left(f(X_{k\delta}^\epsilon,
\hat{Y}_s^\epsilon)-\bar{f}(X_{k\delta}^\epsilon), e_i\right)_H
ds\right]^2\\
&&+2\sum\limits_{i=1}^
\infty\left[\int_{k\delta}^{(k+1)\delta}\sin\{\sqrt{\alpha_i}
s\}\left(f(X_{k\delta}^\epsilon,
\hat{Y}_s^\epsilon)-\bar{f}(X_{k\delta}^\epsilon), e_i\right)_H
ds\right]^2\\
&=&2\sum\limits_{i=1}^
\infty\left[\int_{0}^{\delta}\cos\{\sqrt{\alpha_i}
(s+k\delta)\}\left(f(X_{k\delta}^\epsilon,
\hat{Y}_{s+k\delta}^\epsilon)-\bar{f}(X_{k\delta}^\epsilon),
e_i\right)_H
ds\right]^2\\
&&+2\sum\limits_{i=1}^
\infty\left[\int_{0}^{\delta}\sin\{\sqrt{\alpha_i}
(s+k\delta)\}\left(f(X_{k\delta}^\epsilon,
\hat{Y}_{s+k\delta}^\epsilon)-\bar{f}(X_{k\delta}^\epsilon),
e_i\right)_H ds\right]^2.
\end{eqnarray*}
This means
\begin{eqnarray*}
&&\!\!\!\!\!\!\!\!\mathbb{E}\mathcal{A}_1(t,\epsilon)\\
&\leq&  \sum\limits_{k=0}^{\lfloor
T_0/\delta\rfloor-1}\mathbb{E}\|\mathbbm{i}(t,k,\epsilon)\|_1^2\\
&\leq&{C}\sum\limits_{k=0}^{\lfloor
T_0/\delta\rfloor-1}\sum\limits_{i=1}^
\infty\mathbb{E}\left[\int_{0}^{\delta}\cos\{\sqrt{\alpha_i}
(s+k\delta)\}\left(f(X_{k\delta}^\epsilon,
\hat{Y}_{s+k\delta}^\epsilon)-\bar{f}(X_{k\delta}^\epsilon),
e_i\right)_H ds\right]^2\\
&+&{C}\sum\limits_{k=0}^{\lfloor
T_0/\delta\rfloor-1}\sum\limits_{i=1}^
\infty\mathbb{E}\left[\int_{0}^{\delta}\sin\{\sqrt{\alpha_i}
(s+k\delta)\}\left(f(X_{k\delta}^\epsilon,
\hat{Y}_{s+k\delta}^\epsilon)-\bar{f}(X_{k\delta}^\epsilon),
e_i\right)_H ds\right]^2.
\end{eqnarray*}
Note that, by the construction of $\hat{Y}^\epsilon_t$, for  any
fixed $k$ and $s\in [0, \delta)$  we have

\begin{eqnarray}
\hat{Y}_{s+k\delta}^\epsilon&=&G_{s/\epsilon}{Y}_{k\delta}^\epsilon+\frac{1}{\epsilon}\int_{k\delta}^{k\delta+s}G_{(k\delta+s-r)/\epsilon}g(
X_{k\delta}^\epsilon,\hat{Y}_r^\epsilon)dr\nonumber\\
&&+\frac{1}{\sqrt{\epsilon}}\int_{k\delta}^{k\delta+s}G_{(k\delta+s-r)/\epsilon}b(
X_{k\delta}^\epsilon,\hat{Y}_r^\epsilon)dW_r^{2}.\nonumber
\end{eqnarray}
Taking a time-shift transformation yields
\begin{eqnarray}
\hat{Y}_{s+k\delta}^\epsilon&=&
G_{s/\epsilon}{Y}_{k\delta}^\epsilon+\frac{1}{\epsilon}\int_0^sG_{(s-r)/\epsilon}g(
X_{k\delta}^\epsilon,\hat{Y}_{r+k\delta}^\epsilon) dr\nonumber\\
&&+\frac{1}{\sqrt{\epsilon}}\int_0^sG_{(s-r)/\epsilon}b(
X_{k\delta}^\epsilon,\hat{Y}_{r+k\delta}^\epsilon)dW_r^{*2},\label{Scale-1}
\end{eqnarray}
where $W^{*2}_t$ is the shift version of $W^{2}_t$ and hence they
have the same distribution. Let $\bar{W}_t$ be a $Q_2-$Wiener
process defined on the   stochastic basis  $(\Omega, \mathscr
{F},\mathscr {F}_t,\mathbb{P})$, which is independent of  $W_t^{1}$
and $W_t^{2}$. Denote by $Y^{X_{k\delta}^\epsilon,
{Y}_{k\delta}^\epsilon}$ the unique $H-$valued process satisfying
\begin{eqnarray}
\nonumber Y_{s/\epsilon}^{X_{k\delta}^\epsilon,
{Y}_{k\delta}^\epsilon}&=&G_{s/\epsilon}{Y}_{k\delta}^\epsilon+\int_0^{s/\epsilon}G_{(s/\epsilon-r)}g(
X_{k\delta}^\epsilon,Y_r^{X_{k\delta}^\epsilon,
{Y}_{k\delta}^\epsilon})dr\\
\nonumber&&\qquad\quad\,\,\,+\int_0^{s/\epsilon}G_{(s/\epsilon-r)}b(
X_{k\delta}^\epsilon,Y_r^{X_{k\delta}^\epsilon,
{Y}_{k\delta}^\epsilon})d\bar{W}_r.
\end{eqnarray}
With a simple time-rescaling, we obtain
\begin{eqnarray}
\nonumber Y_{s/\epsilon}^{X_{k\delta}^\epsilon,
{Y}_{k\delta}^\epsilon}&=&
G_{s/\epsilon}{Y}_{k\delta}^\epsilon+\frac{1}{\epsilon}\int_0^sG_{(s-r)/\epsilon}
g(X_{k\delta}^\epsilon, Y^{X_{k\delta}^\epsilon,
{Y}_{k\delta}^\epsilon}_{r/\epsilon})
dr\nonumber\\
&&\qquad\quad\,\,\,+\frac{1}{\sqrt{\epsilon}}\int_0^sG_{(s-r)/\epsilon}b(X_{k\delta}^\epsilon,
Y^{X_{k\delta}^\epsilon,
{Y}_{k\delta}^\epsilon}_{r/\epsilon})d\bar{\bar{W}}_r^\epsilon,\label{Scale-2}
\end{eqnarray}
where $\bar{\bar{W}}_t^\epsilon$ is the scaled version of
$\bar{W}_t$. By comparison, \eqref{Scale-1} and \eqref{Scale-2}
yield
\begin{equation}\label{distribution}
(X^\epsilon_{k\delta},\hat{Y}_{s+k\delta}^\epsilon)\sim
(X^\epsilon_{k\delta}, Y_{s/\epsilon}^{X_{k\delta}^\epsilon,
{Y}_{k\delta}^\epsilon}), \quad s\in [0, \delta),
\end{equation}
where $\sim$ denotes coincidence in distribution sense. In view of
\eqref{distribution} we have
\begin{eqnarray*}
\mathbb{E}\mathcal{A}_1(t,\epsilon)&\leq&{C}
\sum\limits_{k=0}^{\lfloor
T_0/\delta\rfloor-1}(\mathcal{J}_k^\epsilon+\tilde{\mathcal{J}}_k^\epsilon)\\
&\leq& \frac{C}{\delta}\max\limits_{0\leq k\leq \lfloor
T_0/\delta\rfloor}(\mathcal{J}_k^\epsilon+\tilde{\mathcal{J}}_k^\epsilon),
\end{eqnarray*}
here
\begin{eqnarray*}
&&\mathcal {J}_k^\epsilon:=\sum\limits_{i=1}^
\infty\mathbb{E}\left[\int_{0}^{\delta}\sin\{\sqrt{\alpha_i}
(s+k\delta)\}\left(f(X_{k\delta}^\epsilon,
{Y}_{s/\epsilon}^{X_{k\delta}^\epsilon,Y^\epsilon_{k\delta}})-\bar{f}(X_{k\delta}^\epsilon),
e_i\right)_H ds\right]^2,\\
&&\tilde{\mathcal{J}}^\epsilon_k:=\sum\limits_{i=1}^
\infty\mathbb{E}\left[\int_{0}^{\delta}\cos\{\sqrt{\alpha_i}
(s+k\delta)\}\left(f(X_{k\delta}^\epsilon,
{Y}_{s/\epsilon}^{X_{k\delta}^\epsilon,Y^\epsilon_{k\delta}})-\bar{f}(X_{k\delta}^\epsilon),
e_i\right)_H ds\right]^2
\end{eqnarray*}
for $k=0,1,\ldots \lfloor T_0/\delta\rfloor-1.$ In order to prove
Lemma \ref{6-2}, we shall need the following lemma, whose proof can
be founded in Fu et al. \cite[Subsection 6.1]{Fu-Liu-2}
\begin{lemma}\label{6-3}
Suppose that conditions in Lemma \ref{Key-lemma} hold, then there
exists a constant $C>0$ such that
\begin{eqnarray}\label{Claim}
\mathcal {J}_k^\epsilon\leq C\delta\epsilon
\end{eqnarray}
and
\begin{eqnarray}\tilde{\mathcal{J}}^\epsilon_k\leq
C\delta\epsilon\label{Claim-1}
\end{eqnarray}
for $ k=0,1,\cdots, \lfloor T_0/\delta\rfloor-1.$
\end{lemma}

Now, thanks to \eqref{Claim} and \eqref{Claim-1}, we have
\begin{eqnarray}
\mathbb{E}\mathcal{A}_1(t,\epsilon)\leq C {\epsilon}
.\label{I_1-bound}
\end{eqnarray}
Next, by using the idea introduced in Br\'{e}hier \cite{Brehier} and
Liu \cite{LiuDi}, let us estimate
$\mathbb{E}\mathcal{A}_2(t,\epsilon)$. We introduce Markov processes
that generalize $\hat{Y}_t^\epsilon.$ For $k=0,1,\cdots,\lfloor
T_0/\delta\rfloor$, we denote by $\{Z_t^{k,\epsilon}\}_{t\geq
k\delta}$ the solution of the problem
\begin{eqnarray*}
&&{dZ_t^{k,\epsilon}}=\frac{1}{\epsilon}\big(A
Z_t^{k,\epsilon}+g(X_{k\delta}^\epsilon,Z_t^{k,\epsilon})\big)dt\nonumber\label{Z-equation}\\
&&\qquad\quad+\frac{1}{\sqrt{\epsilon}}b(X_{k\delta}^\epsilon,
Z_t^{k,\epsilon})d{W}^{2}_t, \;t\geq k\delta,\\
&&Z_{k\delta}^{k,\epsilon}=\hat{Y}_{k\delta}^\epsilon.
\label{Z-equation-initial}
\end{eqnarray*}
 It is immediate to
check that if $\tau\in[k\delta, (k+1)\delta),$ we have
\begin{equation*}
\hat{Y}^\epsilon_{\tau}=Z^{k,\epsilon}_{\tau}.\label{Z_k-1-coindence}
\end{equation*}
Also, the continuity implies that
\begin{equation}\nonumber
Z^{k,\epsilon}_{(k+1)\delta}=\hat{Y}^\epsilon_{(k+1)\delta}=Z^{k+1,\epsilon}_{(k+1)\delta}.
\end{equation}
For $i\delta\leq s\leq (i+1)\delta\leq j\delta\leq \tau\leq
(j+1)\delta$, we have
\begin{eqnarray*}
&&\Big\langle
\mathbbm{i}(t,i,\epsilon),\mathbbm{i}(t,j,\epsilon)\Big\rangle_1\\
&&=\Big\langle\int_{i\delta}^{(i+1)\delta}S_{t-s}[f(X_{i\delta}^\epsilon,
\hat{Y}_s^\epsilon)-\bar{f}(X_{i\delta}^\epsilon)]ds,\\
&&\qquad\qquad\qquad\qquad\qquad\int_{j\delta}^{(j+1)\delta}S_{t-\tau}[f(X_{j\delta}^\epsilon,
\hat{Y}_\tau^\epsilon)-\bar{f}(X_{j\delta}^\epsilon)]d\tau\Big\rangle_1\\
&&=\sum\limits_{n=1}^\infty\bigg\{\sqrt{\alpha_n}\int_{i\delta}^{(i+1)\delta}\Big(S_{t-s}[f(X_{i\delta}^\epsilon,
\hat{Y}_s^\epsilon)-\bar{f}(X_{i\delta}^\epsilon)],e_n\Big)_Hds\\
&&\qquad\qquad\qquad\qquad\times\sqrt{\alpha_n}\int_{j\delta}^{(j+1)\delta}\Big(S_{t-\tau}[f(X_{j\delta}^\epsilon,
\hat{Y}_\tau^\epsilon)-\bar{f}(X_{j\delta}^\epsilon)],e_n\Big)_Hd\tau\bigg\}\\
&&=\int_{i\delta}^{(i+1)\delta}\int_{j\delta}^{(j+1)\delta}\sum\limits_{n=1}^\infty\bigg[\alpha_n\Big(S_{t-s}[f(X_{i\delta}^\epsilon,
\hat{Y}_s^\epsilon)-\bar{f}(X_{i\delta}^\epsilon)],e_n\Big)_H\\
&&\qquad\qquad\qquad\qquad\qquad\times\Big(S_{t-\tau}[f(X_{j\delta}^\epsilon,
\hat{Y}_\tau^\epsilon)-\bar{f}(X_{j\delta}^\epsilon)],e_n\Big)_H\bigg]dsd\tau.
\end{eqnarray*}
From this, one sees that
\begin{eqnarray*}
&&\mathbb{E}\Big\langle
\mathbbm{i}(t,i,\epsilon),\mathbbm{i}(t,j,\epsilon)\Big\rangle_1\\
&&=\int_{i\delta}^{(i+1)\delta}\int_{j\delta}^{(j+1)\delta}\mathbb{E}\sum\limits_{n=1}^\infty\bigg\{\alpha_n\Big(S_{t-s}[f(X_{i\delta}^\epsilon,
\hat{Y}_s^\epsilon)-\bar{f}(X_{i\delta}^\epsilon)],e_n\Big)_H\\
&&\qquad\qquad\qquad\qquad\times\mathbb{E}\Big[\Big(S_{t-\tau}[f(X_{j\delta}^\epsilon,
\hat{Y}_\tau^\epsilon)-\bar{f}(X_{j\delta}^\epsilon)],e_n\Big)_H\Big|\mathscr{F}_{(i+1)\delta}\Big]\bigg\}dsd\tau\\
&&=\int_{i\delta}^{(i+1)\delta}\int_{j\delta}^{(j+1)\delta}\mathbb{E}\sum\limits_{n=1}^\infty\bigg\{\alpha_n\Big(S_{t-s}[f(X_{i\delta}^\epsilon,
\hat{Y}_s^\epsilon)-\bar{f}(X_{i\delta}^\epsilon)],e_n\Big)_H\\
&&\qquad\qquad\qquad\qquad\times\Big(\mathbb{E}\Big[S_{t-\tau}[f(X_{j\delta}^\epsilon,
\hat{Y}_\tau^\epsilon)-\bar{f}(X_{j\delta}^\epsilon)]\Big|\mathscr{F}_{(i+1)\delta}\Big],e_n\Big)_H\bigg\}dsd\tau\\
&&=\int_{i\delta}^{(i+1)\delta}\int_{j\delta}^{(j+1)\delta}\mathbb{E}\Big\langle
S_{t-s}[f(X_{i\delta}^\epsilon,
\hat{Y}_s^\epsilon)-\bar{f}(X_{i\delta}^\epsilon)],\\
&&\qquad\qquad\qquad\qquad\qquad\mathbb{E}\Big[S_{t-\tau}[f(X_{j\delta}^\epsilon,
\hat{Y}_\tau^\epsilon)-\bar{f}(X_{j\delta}^\epsilon)]\Big|\mathscr{F}_{(i+1)\delta}\Big]\Big\rangle_1dsd\tau\\
&&=\int_{i\delta}^{(i+1)\delta}\int_{j\delta}^{(j+1)\delta}\mathbb{E}\Big\langle
S_{t-s}[f(X_{i\delta}^\epsilon,
\hat{Y}_s^\epsilon)-\bar{f}(X_{i\delta}^\epsilon)],\\
&&\qquad\qquad\qquad\qquad\qquad
S_{t-\tau}\mathbb{E}\Big[f(X_{j\delta}^\epsilon,
\hat{Y}_\tau^\epsilon)-\bar{f}(X_{j\delta}^\epsilon)\Big|\mathscr{F}_{(i+1)\delta}\Big]\Big\rangle_1dsd\tau.
\end{eqnarray*}
Hence, by {boundness} of $f$ we have
\begin{eqnarray*}
&&\mathbb{E}\Big\langle
\mathbbm{i}(t,i,\epsilon),\mathbbm{i}(t,j,\epsilon)\Big\rangle_1\\
&&\leq
C\int_{i\delta}^{(i+1)\delta}\int_{j\delta}^{(j+1)\delta}\mathbb{E}\Big\|S_{t-\tau}\mathbb{E}\Big[f(X_{j\delta}^\epsilon,
\hat{Y}_\tau^\epsilon)-\bar{f}(X_{j\delta}^\epsilon)\Big|\mathscr{F}_{(i+1)\delta}\Big]\Big\|_1dsd\tau\\
&&\leq
C\int_{i\delta}^{(i+1)\delta}\int_{j\delta}^{(j+1)\delta}\mathbb{E}\Big\|\mathbb{E}\Big[f(X_{j\delta}^\epsilon,
\hat{Y}_\tau^\epsilon)-\bar{f}(X_{j\delta}^\epsilon)\Big|\mathscr{F}_{(i+1)\delta}\Big]\Big\|dsd\tau.\\
&&\leq
C\delta\int_{j\delta}^{(j+1)\delta}\mathbb{E}\Big\|\mathbb{E}\Big[f(X_{j\delta}^\epsilon,
\hat{Y}_\tau^\epsilon)-\bar{f}(X_{j\delta}^\epsilon)\Big|\mathscr{F}_{(i+1)\delta}\Big]\Big\|d\tau.
\end{eqnarray*}
Since $\hat{Y}^\epsilon_{\tau}=Z^{j,\epsilon}_{\tau}$, we can
dominate the above integrand  by considering the estimate
\begin{eqnarray}
&&\!\!\!\!\mathbb{E}\Big\|\mathbb{E}\Big[f(X_{j\delta}^\epsilon,
\hat{Y}_\tau^\epsilon)-\bar{f}(X_{j\delta}^\epsilon)\Big|\mathscr{F}_{(i+1)\delta}\Big]\Big\|\nonumber\\
&&\leq\mathbb{E}\Big\|\mathbb{E}\Big[f(X^\epsilon_{(i+1)\delta},Z^{i+1,\epsilon}_\tau)-\bar{f}(X^\epsilon_{(i+1)\delta})
\Big|\mathscr{F}_{(i+1)\delta}]\Big\|\nonumber\\
&&\quad+\mathbb{E}\Big\|\mathbb{E}\Big[\Big(f(X_{j\delta}^\epsilon,
Z^{j,\epsilon}_\tau)-\bar{f}(X_{j\delta}^\epsilon)
\Big)-\Big(f(X_{(i+1)\delta}^\epsilon,
Z^{i+1,\epsilon}_\tau)\nonumber\\
&&\quad-\bar{f}(X_{(i+1)\delta}^\epsilon)
\Big)\Big|\mathscr{F}_{(i+1)\delta}\Big]\Big\|\nonumber\\
&&:=B(i, \tau)+\tilde{B}(i, j,\tau).\nonumber
\end{eqnarray}
Consequently, we have
\begin{eqnarray}
\mathbb{E}\mathcal{A}_2(t,\epsilon)&=&\sum\limits_{0\leq i< j\leq
\lfloor t/\delta\rfloor-1}\mathbb{E}\Big\langle
\mathbbm{i}(t,i,\epsilon),\mathbbm{i}(t,j,\epsilon)\Big\rangle_1\nonumber\\
&\leq&C\delta\sum\limits_{0\leq i< j\leq \lfloor
t/\delta\rfloor-1}\int_{j\delta}^{(j+1)\delta}(B(i,
\tau)+\tilde{B}(i, j,\tau))d\tau.\label{new-1}
\end{eqnarray}
In what follows, we denote by $\{Z^{\epsilon,x,y}_t\}_{t\geq0}$ the
solution of equation
\begin{eqnarray*}
dZ_t^\epsilon=\frac{1}{\epsilon}\big(A
Z_t^\epsilon+g(x,Z_t^\epsilon)\big)dt +\frac{1}{\sqrt{\epsilon}}b(x,
Z_t^\epsilon)d{\tilde{W}}_t,\quad Z_0^\epsilon=y,
\end{eqnarray*}
where the $Q_2-$Wiener processes $\tilde{W}_t$ is independent of
$W^1_t$ and $W^2_t$. It is clear that, for any $k=0,1,\cdots,\lfloor
T_0/\delta\rfloor$, the distribution of the process
$$Z^{k,\epsilon}_t,\quad t\geq k\delta,$$ coincides with the distribution of the
process
$$Z^{\epsilon,X^\epsilon_{k\delta},\hat{Y}^\epsilon_{k\delta}}_{t-k\delta},\quad t\geq k\delta.$$
Hence, thanks to  Markov property, for $(i+1)\delta\leq \tau$ we
obtain
\begin{eqnarray}
B(i,
\tau)&=&\mathbb{E}\Big\|\mathbb{E}\Big[f(x,Z^{\epsilon,x,y}_{{\tau-(i+1)\delta}})-\bar{f}(x)
\Big]\Big|_{x=X_{(i+1)\delta}^\epsilon,y=\hat{Y}_{(i+1)\delta}^\epsilon}\Big\|\nonumber\\
&=&\mathbb{E}\Big\|\mathbb{E}\Big[f(x,Y^{x,y}_{\frac{\tau-(i+1)\delta}{\epsilon}})-\bar{f}(x)
\Big]\Big|_{x=X_{(i+1)\delta}^\epsilon,y=\hat{Y}_{(i+1)\delta}^\epsilon}\Big\|,\nonumber
\end{eqnarray}
where $Y^{x,y}$ denotes the solution of problem
\eqref{Frozen-Fast}-\eqref{Frozen-Initial} with fixed slow component
$x$ and initial value $y$. Then by \eqref{Averaging-Expectation}, we
have
\begin{eqnarray}
B(i, \tau)\leq
Ce^{-c\frac{\tau-(i+1)\delta}{2\epsilon}}(1+\mathbb{E}\|X_{(i+1)\delta}^\epsilon\|+\mathbb{E}\|\hat{Y}_{(i+1)\delta}^\epsilon\|).\nonumber
\end{eqnarray}
Thus, by Lemma \ref{5-1} and Lemma \ref{Y-hat}, for $(i+1)\delta\leq
\tau$ we have
\begin{eqnarray}
B(i, \tau)\leq
Ce^{-c\frac{\tau-(i+1)\delta}{2\epsilon}}.\label{new-2}
\end{eqnarray}
Next, let us estimate $\tilde{B}(i, j,\tau)$. Thanks to the tower
property of conditional expectation, we obtain
\begin{eqnarray}
\tilde{B}(i,
j,\tau)&=&\mathbb{E}\Big\|\mathbb{E}\bigg(\mathbb{E}\Big[\sum\limits_{m=i+1}^{j-1}
[f(X^\epsilon_{(m+1)\delta},Z^{m+1,\epsilon}_\tau)-\bar{f}(X^\epsilon_{(m+1)\delta})]\big|\mathscr{F}_{(m+1)\delta}\Big]\Big|\mathscr{F}_{(i+1)\delta}\bigg)\nonumber\\
&&\quad-\mathbb{E}\bigg(\mathbb{E}\Big[\sum\limits_{m=i+1}^{j-1}
[f(X^\epsilon_{m\delta},Z^{m,\epsilon}_\tau)-\bar{f}(X^\epsilon_{m\delta})]\big|\mathscr{F}_{(m+1)\delta}\Big]\Big|\mathscr{F}_{(i+1)\delta}\bigg)
\Big\|\nonumber\\
&\leq&\sum\limits_{m=i+1}^{j-1}\mathbb{E}\Big\|\mathbb{E}\Big[
f(X^\epsilon_{(m+1)\delta},Z^{m+1,\epsilon}_\tau)-\bar{f}(X^\epsilon_{(m+1)\delta})\big|\mathscr{F}_{(m+1)\delta}\Big]\nonumber\\
&&\qquad\qquad-\mathbb{E}\Big[
f(X^\epsilon_{m\delta},Z^{m,\epsilon}_\tau)-\bar{f}(X^\epsilon_{m\delta})\big|\mathscr{F}_{(m+1)\delta}\Big]\Big\|.\nonumber
\end{eqnarray}
Due to the Markov property, we have
\begin{eqnarray}
&&\mathbb{E}\Big[
f(X^\epsilon_{(m+1)\delta},Z^{m+1,\epsilon}_\tau)-\bar{f}(X^\epsilon_{(m+1)\delta})
\big| \mathscr{F}_{(m+1)\delta}\Big]\nonumber\\
&&=\mathbb{E}\Big[f(x,Z^{\epsilon,x,y}_{{\tau-(m+1)\delta}})-\bar{f}(x)
\Big]{\Big|_{x=X_{(m+1)\delta}^\epsilon,y=\hat{Y}_{(m+1)\delta}^\epsilon}}\nonumber\\
&&=\mathbb{E}\Big[f(x,Y^{x,y}_\frac{{\tau-(m+1)\delta}}{\epsilon})-\bar{f}(x)
\Big]{\Big|_{x=X_{(m+1)\delta}^\epsilon,y=\hat{Y}_{(m+1)\delta}^\epsilon}}.\nonumber
\end{eqnarray}
Similarly, one has
\begin{eqnarray}
&&\mathbb{E}\Big[
f(X^\epsilon_{m\delta},Z^{m,\epsilon}_\tau)-\bar{f}(X^\epsilon_{m\delta})
\big| \mathscr{F}_{(m+1)\delta}\Big]\nonumber\\
&&=\mathbb{E}\Big[f(x,Z^{\epsilon,x,y}_{{\tau-(m+1)\delta}})-\bar{f}(x)
\Big]{\Big|_{x=X_{m\delta}^\epsilon,y=\hat{Y}_{(m+1)\delta}^\epsilon}}\nonumber\\
&&=\mathbb{E}\Big[f(x,Y^{x,y}_\frac{{\tau-(m+1)\delta}}{\epsilon})-\bar{f}(x)
\Big]{\Big|_{x=X_{m\delta}^\epsilon,y=\hat{Y}_{(m+1)\delta}^\epsilon}}.\nonumber
\end{eqnarray}
Let us define
\begin{eqnarray}
\tilde{f}(x,y,t):=\mathbb{E}f(x, Y_t^{x,y})-\bar{f}(x), \; x,y \in
H.\label{tilde-f}
\end{eqnarray}
 Hence, thanks to
Lemma \ref{Lemma7.1} presented in the final section, for
$(i+1)\delta\leq j\delta\leq \tau\leq (j+1)\delta$ we have
\begin{eqnarray}
\tilde{B}(i, j,\tau)&\leq&\sum\limits_{m=i+1}^{j-1}\mathbb{E}\Big\|\tilde{f}(X_{m\delta}^\epsilon,\hat{Y}_{(m+1)\delta}^\epsilon,\frac{{\tau-(m+1)\delta}}{\epsilon})\nonumber\\
&&\qquad\qquad-\tilde{f}(X_{(m+1)\delta}^\epsilon,\hat{Y}_{(m+1)\delta}^\epsilon,\frac{{\tau-(m+1)\delta}}{\epsilon})\Big\|\nonumber\\
&&\leq
C\sum\limits_{m=i+1}^{j-1}e^{-c\frac{{\tau-(m+1)\delta}}{\epsilon}}
\big(\mathbb{E}\|X_{(m+1)\delta}^\epsilon-X_{m\delta}^\epsilon\|\big)\big(1+\mathbb{E}\|\hat{Y}_{(m+1)\delta}^\epsilon\|\big)\nonumber\\
&&\leq
C\delta\sum\limits_{m=i+1}^{j-1}e^{-c\frac{{\tau-(m+1)\delta}}{\epsilon}}\nonumber\\
&&\leq C\delta
\frac{e^{-c\frac{\tau-j\delta}{\epsilon}}}{1-e^{-c\frac{\delta}{\epsilon}}},\label{new-3}
\end{eqnarray}
where we used \eqref{Regularity} and \eqref{Y-hat-1} in the second
to last step. Now, recalling \eqref{new-1}, due to  \eqref{new-2}
and \eqref{new-3},  we get
\begin{eqnarray}
\mathbb{E}\mathcal{A}_2(t,\epsilon)&\leq& C\sum\limits_{0\leq i<
j\leq \lfloor
t/\delta\rfloor-1}\delta\int_{j\delta}^{(j+1)\delta}e^{-c\frac{\tau-(i+1)\delta}{\epsilon}}d\tau\nonumber\\
&&+C\sum\limits_{0\leq i< j\leq \lfloor
t/\delta\rfloor-1}\delta\int_{j\delta}^{(j+1)\delta}\delta
\frac{e^{-c\frac{\tau-j\delta}{\epsilon}}}{1-e^{-c\frac{\delta}{\epsilon}}}d\tau\nonumber\\
&&\leq C\sum\limits_{0\leq i< j\leq \lfloor
t/\delta\rfloor-1}\left[e^{-\frac{c}{\epsilon}(j-i-1)\delta}(1-e^{-\frac{c}{\epsilon}\delta})\delta\epsilon+\delta^2\epsilon\right]\nonumber\\
&&\leq C\epsilon.\nonumber
\end{eqnarray}
This fact, together with estimate \eqref{I_1-bound} and equality
\eqref{I_1-inner}, shows
\begin{eqnarray}
\mathbb{E}\|I(t,\epsilon)\|^2_1\leq C\epsilon.
\end{eqnarray}
\textbf{Step 4:} Estimate of $
\mathbb{E}\|{\mathcal{L}_t^\epsilon}\|^2_1$.

It is now easy to gather all previous estimates for terms in
\eqref{I_1--I_3} and deduce
\begin{eqnarray}
\mathbb{E}\|\mathcal {L}^\epsilon_t\|^2_1\leq
C(\delta^2+{\epsilon}).\nonumber
\end{eqnarray}
\textbf{Step 5:} Estimate of $
\mathbb{E}\|\dot{\mathcal{L}_t^\epsilon}\|^2.$

It is also easy to see that
\begin{eqnarray}
\dot{\mathcal{L}_t^\epsilon}&=&\sum\limits_{i=1}^\infty
e_i\cdot\int_0^t\cos\{\sqrt{\alpha_i}(t-s)\}\Big(f(X_{s(\delta)}^\epsilon,
\hat{Y}_{s}^\epsilon)-\bar{f}(X_s^\epsilon),e_i\Big)_Hds\nonumber\\
&=&\int_0^tS_{t-s}'\left(f(X_{s(\delta)}^\epsilon,
\hat{Y}_{s}^\epsilon)-\bar{f}(X_s^\epsilon)\right)ds,\nonumber
\end{eqnarray}
and so that we can   decompose $\dot{\mathcal{L}_t^\epsilon}$ as
\begin{equation}\label{J_1--J_3}
\dot{\mathcal {L}}^\epsilon_t= \tilde{I}_1({t},
\epsilon)+\tilde{I}_2({t}, \epsilon)+\tilde{I}_3({t}, \epsilon),
\end{equation}
where
\begin{eqnarray*}
&&\tilde{I}_1({t}, \epsilon)=\sum\limits_{k=0}^{\lfloor
t/\delta\rfloor-1}\int_{k\delta}^{(k+1)\delta}
S_{t-s}'\left(f(X_{k\delta}^\epsilon,
\hat{Y}_s^\epsilon)-\bar{f}(X_{k\delta}^\epsilon)\right)ds,\\
&&\tilde{I}_2({t}, \epsilon)=\sum\limits_{k=0}^{\lfloor
t/\delta\rfloor-1}\int_{k\delta}^{(k+1)\delta}
S_{t-s}'\left(\bar{f}(X_{k\delta}^\epsilon)
-\bar{f}(X_s^\epsilon)\right)ds\\
&&\qquad\quad=\int_0^{\lfloor
t/\delta\rfloor\delta}S_{t-s}' \left(\bar{f}(X_{s(\delta)}^\epsilon)-\bar{f}(X_s^\epsilon)\right)ds,\\
&&\tilde{I}_3({t}, \epsilon)=\int_{{\lfloor
t/\delta\rfloor}\delta}^{t} S_{t-s}'\left(f(X_{{\lfloor
t/\delta\rfloor}\delta}^\epsilon,
\hat{Y}^\epsilon_s)-\bar{f}(X_s^\epsilon)\right)ds.
\end{eqnarray*}
For $\tilde{I}_2({t}, \epsilon)$, using the H\"{o}lder inequality
and the Lipschitz continuity of $\bar{f}$ yields
\begin{eqnarray}
\|\tilde{I}_2(t, \epsilon)\|^2&\leq& \left[\int_0^{{\lfloor
t/\delta\rfloor}\delta}\left\|S_{t-s}'\left(\bar{f}(X_{s(\delta)}^\epsilon)-\bar{f}(X_s^\epsilon)\right)\right\|ds\right]^2\nonumber\\
&\leq&T_0\int_0^{T_0}\left\|\bar{f}(X_{s(\delta)}^\epsilon)-\bar{f}(X_s^\epsilon)\right\|^2ds\nonumber\\
&\leq&C\int_0^{T_0}\left\|X_{s(\delta)}^\epsilon-X_s^\epsilon\right\|^2ds.\nonumber
\end{eqnarray}
By virtue of \eqref{Regularity} in Lemma \ref{5-2} we obtain
\begin{equation}
\mathbb{E}\|\tilde{I}_2(t, \epsilon)\|^2\leq
C\delta^2.\label{J_2-bound}
\end{equation}
Concerning $\tilde{I}_3({t}, \epsilon)$,   by means of  the
H\"{o}lder inequality we get
\begin{eqnarray}
\|\tilde{I}_3({t}, \epsilon)\|^2&\leq& \delta \int_{{\lfloor
t/\delta\rfloor}\delta}^{t} \left\|S_{t-s}'\left(f(X_{{\lfloor
t/\delta\rfloor}\delta}^\epsilon,
\hat{Y}^\epsilon_s)-\bar{f}(X_s^\epsilon)\right)\right\|^2ds \nonumber\\
&\leq&\delta C\int_{\lfloor
t/\delta\rfloor\delta}^{t}\left(1+\|X_s^\epsilon\|^2+\|X_{{\lfloor
t/\delta\rfloor}\delta}^\epsilon\|^2+\|\hat{Y}^\epsilon_s\|^2\right)ds.\nonumber
\end{eqnarray}
Due to the fact that $\|\cdot\|-$norm is bounded by
$\|\cdot\|_1-$norm, we conclude from \eqref{Boundness-slow-solution}
and   \eqref{Y-hat-1} that
\begin{eqnarray}
\mathbb{E}\|\tilde{I}_3({t}, \epsilon)\|^2&\leq&
C\delta^2.\label{J_3-bound}
\end{eqnarray}
For $\tilde{I}_1({t}, \epsilon)$, we have by the series
representation of the Green's function that
\begin{eqnarray*}
\tilde{I}_1(t, \epsilon)=\sum\limits_{k=0}^{\lfloor
t/\delta\rfloor-1}\tilde{\mathbbm{i}}(t,k,\epsilon),
\end{eqnarray*}
where
\begin{eqnarray*}
&&\!\!\!\!\!\!\!\!\!\!\!\!\tilde{\mathbbm{i}}(t,k,\epsilon)\\
&:=&\int_{k\delta}^{(k+1)\delta}\sum\limits_{i=1}^
\infty{\cos\{\sqrt{\alpha_i}(t-s)\}}\left(f(X_{k\delta}^\epsilon,
\hat{Y}_s^\epsilon)-\bar{f}(X_{k\delta}^\epsilon), e_i\right)_H\cdot
e_i ds\\
&=&\sum\limits_{i=1}^ \infty{\cos\{{\sqrt{\alpha_i}t}\}}\cdot
e_i\cdot \int_{k\delta}^{(k+1)\delta}\cos\{\sqrt{\alpha_i}
s\}\left(f(X_{k\delta}^\epsilon,
\hat{Y}_s^\epsilon)-\bar{f}(X_{k\delta}^\epsilon), e_i\right)_H
ds\\
&+&\sum\limits_{i=1}^ \infty{\sin{\{\sqrt{\alpha_i}t\}}}\cdot
e_i\cdot \int_{k\delta}^{(k+1)\delta}\sin\{\sqrt{\alpha_i}
s\}\left(f(X_{k\delta}^\epsilon,
\hat{Y}_s^\epsilon)-\bar{f}(X_{k\delta}^\epsilon), e_i\right)_H ds
\end{eqnarray*}
for $k=0,1,\cdots,\lfloor t/\delta\rfloor-1.$ It is clear that
\begin{eqnarray}
\|\tilde{I}_1(t, \epsilon)\|^2&=&\sum\limits_{k=0}^{\lfloor
t/\delta\rfloor-1}\|\tilde{\mathbbm{i}}(t,k,\epsilon)\|^2+2\sum\limits_{0\leq
i< j\leq \lfloor t/\delta\rfloor-1}\Big(
\tilde{\mathbbm{i}}(t,i,\epsilon),\tilde{\mathbbm{i}}(t,j,\epsilon)\Big)_H\nonumber\\
&:=&\tilde{\mathcal{A}}_1(t,\epsilon)+2\tilde{\mathcal{A}}_2(t,\epsilon).\nonumber
\end{eqnarray}
Applying similar arguments as in  Step 3, we can conclude that
\begin{eqnarray*}
\mathbb{E}\tilde{\mathcal{A}}_1(t,\epsilon)\leq C {\epsilon}
\end{eqnarray*}
and
\begin{eqnarray*}
\mathbb{E}\tilde{\mathcal{A}}_2(t,\epsilon)\leq C \epsilon,
\end{eqnarray*}
which means
\begin{eqnarray}
\mathbb{E}\|\tilde{I}_1(t, \epsilon)\|^2\leq C
\epsilon\label{J_1-bound}.
\end{eqnarray}
Now, in view of estimates \eqref{J_2-bound}, \eqref{J_3-bound} and
\eqref{J_1-bound}, from \eqref{J_1--J_3} we obtain
\begin{eqnarray*}
\mathbb{E}\|\dot{\mathcal {L}}^\epsilon_t\|^2\leq
C({\epsilon}+\delta^2).
\end{eqnarray*}
\textbf{Step 6:} Conclusion. By making use of the results in Step 4
and Step 5, we get the desired inequality \eqref{sup-A-1}.
\end{proof}
\end{lemma}
\subsection{Proof of Theorem \ref{Theorem}}
Now we are in position to conclude the proof for the main theorem
with the aid of the above lemmas.
\begin{proof}
According to \eqref{Auxiliary-slow-difference},
\eqref{Auxiliary-slow-derive} and  Lemma \ref{Key-lemma}, we have
\begin{eqnarray*}
\mathbb{E}\left
(\|{X}_t^\epsilon-\bar{X}_t\|^{2}_1+\|\dot{{X}}_t^\epsilon-\dot{\bar{X}}_t\|^{2}\right
)\leq C({\epsilon+\delta^2}).
\end{eqnarray*}
Taking $\delta=\sqrt{\epsilon}$, we obtain
\begin{eqnarray*}
\mathbb{E}\left
(\|{X}_t^\epsilon-\bar{X}_t\|^{2}_1+\|\dot{{X}}_t^\epsilon-\dot{\bar{X}}_t\|^{2}\right
)\leq C\epsilon,
\end{eqnarray*}
which completes the proof.
\end{proof}

\section{Auxiliary Lemma}
In this section, we state and prove a technical lemma  used in the
former section.
\begin{lemma}\label{Lemma7.1}
Function $\tilde{f}$ defined by \eqref{tilde-f} is Lipschitz
continuous with respect to $x$. In additional, there exist $c, C>0
$, such that for any $ x_1, x_2, y \in H$ and $t>0$ we have
\begin{eqnarray*}
\|\tilde{f}(x_1,y, t)-\tilde{f}(x_2,y, t)\|\leq
C(1+\|y\|)\|x_1-x_2\|e^{-ct}.
\end{eqnarray*}
\begin{proof}
We shall follow the approach of Br\'{e}hier \cite[Proposition
C.2]{Brehier}. For any $t_0>0$, we set
\begin{eqnarray*}
\tilde{F}_{t_0}(x,y,t)=F(x,y,t)-F(x,y,t+t_0),
\end{eqnarray*}
where
\begin{eqnarray*}
F(x,y,t):=\mathbb{E}f(x, Y^{x,y}_t).
\end{eqnarray*}
Thanks to Markov property we then write that
\begin{eqnarray*}
\tilde{F}_{t_0}(x,y,t)&=&F(x,y,t)-\mathbb{E}f(x,Y_{t+t_0}^{x,y})\\
&=&F(x,y,t)-\mathbb{E}F(x, Y_{t_0}^{x,y},t)
\end{eqnarray*}
In view of the assumption (A1), $F$ is G\^{a}teaux-differentiable
with respect to $x$ at $(x,y,t)$. Therefore, we have for any $h\in
H$ that
\begin{eqnarray}
D_x\tilde{F}_{t_0}(x,y,t)\cdot h&=&D_xF(x,y,t)\cdot
h-\mathbb{E}D_x\left(F(x, Y_{t_0}^{x,y},t)\right)\cdot h\nonumber\\
&=&F_x'(x,y,t)\cdot h-\mathbb{E}F_x'(x, Y_{t_0}^{x,y},t)\cdot
h\nonumber\\
&&-\mathbb{E}F_y'(x,
Y_{t_0}^{x,y},t)\cdot\left(D_xY_{t_0}^{x,y}\cdot
h\right),\label{7-1-1}
\end{eqnarray}
where we use the symbol $F_x'$ and $F_y'$ to denote the G\^{a}teaux
derivative with respect to $x$ and $y$, respectively. Note that the
first derivative $\zeta_t^{x,y, h}=D_xY_{t}^{x,y}\cdot h$, at the
point $x$ and along the direction $h\in H$, is the solution to
equation
\begin{eqnarray*}
d\zeta_t^{x, y,h}&=&\left(A\zeta_t^{x, y,h}+g_x'(x, Y_t^{x,y})\cdot
h+g_y'(x, Y_t^{x,y})\cdot\zeta_t^{x, y,h}\right)dt\\
&&+\left(b_{x}'(x, Y_t^{x,y})\cdot h+b_{y}'(x,
Y_t^{x,y})\cdot\zeta_t^{x, y,h}\right)dW^2_t
\end{eqnarray*}
with initial value $\zeta_0^{x,y, h}=0$. We have
\begin{eqnarray}
\frac{d}{dt}\mathbb{E}\|\zeta_t^{x, y,h}\|^2&=&2\mathbb{E}\langle A
\zeta_t^{x, y,h}, \zeta_t^{x, y,h}\rangle+2\mathbb{E}\big(g_x'(x,
Y_t^{x,y})\cdot h, \zeta_t^{x, y,h}\big)_H\nonumber\\
&+&2\mathbb{E}\big(g_y'(x, Y_t^{x,y})\cdot\zeta_t^{x, y,h},
\zeta_t^{x, y,h}\big)_H\nonumber\\
&+&\mathbb{E}\|b_{x}'(x, Y_t^{x,y})\cdot h+b_{y}'(x,
Y_t^{x,y})\cdot\zeta_t^{x, y,h}\|^2_{Q_2}.\nonumber
\end{eqnarray}
Due to the Poincar\'{e} inequality, we get
\begin{eqnarray}
2 \langle A \zeta_t^{x, y,h}, \zeta_t^{x, y,h}\rangle\leq
-2\alpha_1\|\zeta_t^{x, y,h}\|^2.\label{appen-new-5}
\end{eqnarray}
By the Young inequality, for any $\rho>0$, there exists a constant
$C_\rho$ such that
\begin{eqnarray*}
2\mathbb{E}\Big|\big(g_x'(x, Y_t^{x,y})\cdot h, \zeta_t^{x,
y,h}\big)_H\Big|\leq\rho \mathbb{E}\|\zeta_t^{x,
y,h}\|^2+C_\rho\mathbb{E}\|g_x'(x, Y_t^{x,y})\cdot h\|^2.
\end{eqnarray*}
 {According to condition \eqref{g-1}}, this means
\begin{eqnarray}
2\mathbb{E}\Big|\big(g_x'(x, Y_t^{x,y})\cdot h, \zeta_t^{x,
y,h}\big)_H\Big|\leq\rho \mathbb{E}\|\zeta_t^{x,
y,h}\|^2+C_\rho\|h\|^2.\label{appen-new-6}
\end{eqnarray}
 {By using condition \eqref{g-2}}, we easily obtain
\begin{eqnarray}
2\mathbb{E}\big(g_y'(x, Y_t^{x,y})\cdot\zeta_t^{x, y,h}, \zeta_t^{x,
y,h}\big)_H\leq 2L_g\mathbb{E}\|\zeta_t^{x,
y,h}\|^2.\label{appen-new-7}
\end{eqnarray}
By the Young inequality and conditions \eqref{b-1} and \eqref{b-2},
we have
\begin{eqnarray}
&&\mathbb{E}\|b_{x}'(x, Y_t^{x,y})\cdot h+b_{y}'(x,
Y_t^{x,y})\cdot\zeta_t^{x, y,h}\|^2_{Q_2}\nonumber\\ &&\leq
(1+\rho)L_b^2\mathbb{E}\|\zeta_t^{x, y,h}\|^2
+(1+C_\rho\|h\|^2).\label{appen-new-8}
\end{eqnarray}
Thanks to inequalities \eqref{appen-new-5}, \eqref{appen-new-6},
\eqref{appen-new-7} and \eqref{appen-new-8}, by choosing $\rho$
small enough it follows from condition \eqref{decay} that
\begin{eqnarray*}
\frac{d}{dt}\mathbb{E}\|\zeta_t^{x, y,h}\|^2\leq
-C_1\mathbb{E}\|\zeta_t^{x, y,h}\|^2+C_2\|h\|^2
\end{eqnarray*}
for some constants $C_1, C_2>0$. Hence, thanks to  Gronwall's
inequality, it is immediate to check that for any $t\geq 0$,
\begin{eqnarray}
\mathbb{E}\|\zeta_t^{x,y, h}\|^2\leq C\|h\|^2.\label{7-2}
\end{eqnarray}
Notice that for any $y_1, y_2\in H$, we have
\begin{eqnarray}
\|F(x,y_1, t)-F(x,y_2,t)\|&=&\|\mathbb{E}f(x,
Y_t^{x,y_1})-\mathbb{E}f(x,
Y_t^{x,y_2})\|\nonumber\\
&\leq&C\mathbb{E}\|Y_t^{x,y_1}-Y_t^{x,y_2}\|\nonumber\\
&\leq& Ce^{-ct}\|y_1-y_2\|,\nonumber
\end{eqnarray}
which implies
\begin{eqnarray}
\|F_y'(x, y,t)\cdot k\|\leq Ce^{-ct}\|k\|,\; k\in H.\label{7-3}
\end{eqnarray}
Therefore, thanks to \eqref{7-2} and \eqref{7-3}, we can conclude
that
\begin{eqnarray}
\|\mathbb{E}[F_y'(x,
Y_{t_0}^{x,y},t)\cdot\left(D_xY_{t_0}^{x,y}\cdot h\right)]\|\leq C
e^{-ct}\|h\|.\label{7-4}
\end{eqnarray}
Then, we directly have
\begin{eqnarray}
&&F_x'(x,y_1,t)\cdot h-F_x'(x,y_2,t)\cdot h\nonumber\\
&&\quad=\mathbb{E}\left(f_x'(x, Y_t^{x,y_1})\right)\cdot
h-\mathbb{E}\left(f_x'(x, Y_t^{x,y_2})\right)\cdot h\nonumber\\
&&\quad\quad+\mathbb{E}\left(f_y'(x, Y_t^{x,y_1})\cdot
\zeta_t^{x,y_1,
h}-f_y'(x, Y_t^{x,y_2})\cdot \zeta_t^{x,y_2, h}\right)\nonumber\\
&&\quad= \mathbb{E}\left(f_x'(x, Y_t^{x,y_1})\right)\cdot
h-\mathbb{E}\left(f_x'(x, Y_t^{x,y_2})\right)\cdot h\nonumber\\
&&\quad\quad+\mathbb{E}\left([f_y'(x, Y_t^{x,y_1})-f_y'(x,
Y_t^{x,y_2})]\cdot\zeta_t^{x,y_1, h} \right)\nonumber\\
&&\quad\quad+\mathbb{E}\left(f_y'(x, Y_t^{x,
y_2})\cdot(\zeta_t^{x,y_1, h}-\zeta_t^{x,y_2, h})\right).\label{7-5}
\end{eqnarray}
Firstly, in view of  \eqref{f-3} it is easy to show
\begin{eqnarray}
&&\|\mathbb{E}\left(f_x'(x, Y_t^{x,y_1})\right)\cdot
h-\mathbb{E}\left(f_x'(x, Y_t^{x,y_2})\right)\cdot h\|\nonumber\\
&&\quad\leq\mathbb{E}\|\left(f_x'(x, Y_t^{x,y_1})\right)\cdot
h-\left(f_x'(x, Y_t^{x,y_2})\right)\cdot h\|\nonumber\\
&&\quad\leq
C\mathbb{E}\|Y_t^{x,y_1}-Y_t^{x,y_2}\|\cdot\|h\|\nonumber\\
&&\quad\leq Ce^{-ct}\|y_1-y_2\|\cdot\|h\|.\label{7-6}
\end{eqnarray}
Next, by condition \eqref{f-4}  we have
\begin{eqnarray}
&&\|\mathbb{E}\left([f_y'(x, Y_t^{x,y_1})-f_y'(x,
Y_t^{x,y_2})]\cdot\zeta_t^{x,y_1, h} \right)\|\nonumber\\
&&\quad\leq\mathbb{E}\|[f_y'(x, Y_t^{x,y_1})-f_y'(x,
Y_t^{x,y_2})]\cdot\zeta_t^{x,y_1, h}\|\nonumber\\
&&\quad\leq C\{\mathbb{E}\|\zeta_t^{x,y_1,
h}\|^2\}^{\frac{1}{2}}\cdot\{\mathbb{E}\|Y_t^{x,y_1}-Y_t^{x,y_2}\|^2\}^{\frac{1}{2}}\nonumber\\
&&\quad\leq C e^{-ct}\|h\|\cdot\|y_1-y_2\|.\label{7-7}
\end{eqnarray}
By making use of  condition \eqref{f-2}, we can show
\begin{eqnarray}
&&\|\mathbb{E}\left(f_y'(x, Y_t^{x, y_2})\cdot(\zeta_t^{x,y_1,
h}-\zeta_t^{x,y_2, h})\right)\|\nonumber\\
&&\quad\leq\mathbb{E}\|\left(f_y'(x, Y_t^{x,
y_2})\cdot(\zeta_t^{x,y_1, h}-\zeta_t^{x,y_2, h})\right)\|\nonumber\\
&&\quad\leq C\mathbb{E}\|\zeta_t^{x,y_1, h}-\zeta_t^{x,y_2, h}\|\nonumber\\
&&\quad\leq C e^{-c't}\|y_1-y_2\|\cdot\|h\|.\label{7-8}
\end{eqnarray}
Collecting together \eqref{7-5}, \eqref{7-6}, \eqref{7-7} and
\eqref{7-8}, we get
\begin{eqnarray*}
&&\|F_x'(x,y_1,t)\cdot h-F_x'(x,y_2,t)\cdot h\|\nonumber\\
&&\leq C e^{-c_0t}\|y_1-y_2\|\cdot\|h\|,
\end{eqnarray*}
which means
\begin{eqnarray}
&&\|F_x'(x,y,t)\cdot h-\mathbb{E}F_x'(x,Y^{x,y}_{t_0},t)\cdot h\|\nonumber\\
&&\leq C e^{-c_0t}(1+\|y\|)\cdot\|h\| \label{7-9}
\end{eqnarray}
since
\begin{eqnarray*}
\mathbb{E}\|Y^{x,y}_{t_0}\|\leq C(1+\|y\|).
\end{eqnarray*}
Returning to \eqref{7-1-1}, by \eqref{7-4} and \eqref{7-9} we
conclude that
\begin{eqnarray*}
\|D_x\tilde{F}_{t_0}(x,y,t)\cdot h\|\leq Ce^{-ct}(1+\|y\|)\|h\|.
\end{eqnarray*}
This yields
\begin{eqnarray*}
\|\tilde{F}_{t_0}(x_1,y,t)-\tilde{F}_{t_0}(x_2,y,t)\|\leq
Ce^{-ct}(1+\|y\|)\|x_1-x_2\|.
\end{eqnarray*}
Letting $t_0\rightarrow +\infty$, we obtain
\begin{eqnarray*}
\|\tilde{f}(x_1,y, t)-\tilde{f}(x_2,y, t)\|\leq
C(1+\|y\|)\|x_1-x_2\|e^{-ct}.
\end{eqnarray*}
\end{proof}
\end{lemma}



\end{document}